\documentclass{amsart}
\usepackage{amsmath, amssymb, amsthm}

\newcommand{\llq}{\lq\lq} 

\newcommand{\Z}{\mathbb Z} 
\newcommand{\R}{\mathbb R} 

\newcommand{\dee}{\mathrm d} 
\providecommand{\diff}[2]{\frac{\dee #1}{\dee #2}} 
\providecommand{\diffn}[3]{\frac{\dee^{#1} #2}{\dee #3^{#1}}} 
\newcommand{\del}{\partial} 
\providecommand{\pdiff}[2]{\frac{\del #1}{\del #2}} 
\providecommand{\pdiffn}[3]{\frac{\del^{#1} #2}{\del #3^{#1}}} 

\providecommand{\vect}[1]{\mathbf{#1}} 
\providecommand{\matr}[1]{\mathrm{#1}} 

\newtheorem{thm}{Theorem}[section]
\newtheorem{defn}[thm]{Definition}
\newtheorem{prop}[thm]{Proposition}
\newtheorem{lemma}[thm]{Lemma}
\newtheorem{cor}[thm]{Corollary}
\newtheorem{rmrk}[thm]{Remark}

\numberwithin{equation}{section}
\newcommand{\prf}{\textbf{\underline{Proof}:} }

\begin{document}

\title[Geometry of star-shaped curves]{\bf On the geometry of star-shaped curves in $\R^n$}
\author{Stefan A. Horocholyn}
\address{Present Address (Visiting): Department of Mathematical Sciences, Tokyo Metropolitan University, 1-1 Minami-Osawa, Hachioji-shi, Tokyo, 192-0397, Japan}
\email{shoro@tmu.ac.jp}

\begin{abstract}
The manifold $\mathcal{M}$ of star-shaped curves in $\R^n$ is considered via the theory of connections on vector bundles, and cyclic $\mathcal{D}$-modules. 
The appropriate notion of an $\llq$integral curve" (i.e. certain admissible deformations) on $\mathcal{M}$ is defined, and the resulting space of admissible deformations is classified via iso-spectral flows, which are shown to be described by equations from the $n$-KdV hierarchy. 
\end{abstract}


\maketitle

\section{Introduction}

Smooth, 1-parameter deformations of curves, and the induced evolution equations satisfied by the differential invariants of the curve (curvature, torsion, etc.) are well-known to be closely related to the theory of solitons (e.g. \cite{Lamb, Hasimoto, CQ, Ino10, Musso, Pinkall, Terng97}). In this article, we characterize such $\llq$integrable deformations" for the class of star-shaped curves in $\R^n$ (defined below), by showing that the relevant curvature functions must satisfy an equation of the $n$-KdV hierarchy. This builds on results for the cases $n=2$ and $n=3$ for equi-centro-affine curves \cite{BuPePi, CIMB09, CIMB13, CQ, FujiKuro10, FujiKuro14, GLHMB, HuangSinger, Ino10, Pinkall, TerngWu13}; some of our results parallel the preprint \cite{TerngWu15}, which builds on work by Terng-Uhlenbeck on applications of loop groups to integrable systems (see \cite{TU00, TU11} and references therein). 

Despite the different approaches considered so far in the literature, the core object is a flat connection on a smooth vector bundle, with a distinguished section identified as a framed curve, and its deformations which preserve the rank of the space of smooth sections on the vector bundle. We explain this point-of-view, and in doing so, provide a differential geometric realization of aspects of the theory of the $n$-KdV hierarchy (see e.g. \cite{GelDik, Adler, Dub, DrinSok}).

Let us consider a smooth curve $\gamma:I \rightarrow \R^{1 \times n}$ given by $\gamma = (\gamma_1, \dots, \gamma_n)$, where $I = (a,b)$ is a real interval (possibly infinite), with smooth coordinate $x$ (arbitrary but fixed), and $n > 1$. Letting $\gamma^{(k)} := \diffn{k}{\gamma}{x}$ denote the $k^\matr{th}$ component-wise derivative of $\gamma$, we say that $\gamma$ is \textit{star-shaped} whenever the following matrix is smooth and invertible for all $x \in I$:
\[ \matr{W}(\gamma_1, \dots, \gamma_n) := \begin{bmatrix} \gamma \\ \gamma^{(1)} \\ \vdots \\ \gamma^{(n-1)} \end{bmatrix} 
= \begin{bmatrix} y_1 & y_2 & \dots & y_n \end{bmatrix} \ , \quad 
y_j := \begin{bmatrix} \gamma_j \\ \gamma_j^{(1)} \\ \vdots \\ \gamma_j^{(n-1)} \end{bmatrix} \ . \]
Letting $\mathcal{F} := \mathcal{C}^{\infty}(I, \R)$, this implies that\footnote{Since we are fixing the smooth coordinate $x: I \rightarrow \R$, the notation $\mathcal{F}^n$ is equivalent to the more precise but cumbersome $\mathcal{C}^{\infty}(I, \R^n)$.} $\mathcal{F}^n \cong < \gamma, \gamma^{(1)}, \dots, \gamma^{(n-1)} >_\mathcal{F}$, so $\gamma^{(n)}$ is necessarily a linear combination of these basis vectors, and writing this out in terms of $\matr{W}(\gamma)$ yields:
\begin{equation}\label{FrenetSerret}
\diff{}{x}\matr{W}(\gamma) = \matr{F}\matr{W}(\gamma) \ , \quad
\matr{F} := \begin{bmatrix}   0  & \cdot  & \cdots         & \cdot    \\
                              0  & \vdots & \matr{I}_{n-1} & \vdots   \\
                              0  & \cdot  & \cdots         & \cdot    \\
                            -u_0 &  -u_1  & \cdots         & -u_{n-1}  \end{bmatrix} \ ,
\end{equation}
for some functions $u_k \in \mathcal{F}$. Evidently, $\matr{W}(\gamma)$ is a fundamental solution matrix to the linear homogeneous ODE $L_{\gamma}(y)=0$ of order $n$, where:
\[ L_{\gamma}(y) := \frac{\det\matr{W}(\gamma_1, \dots, \gamma_n, y)}{\det\matr{W}(\gamma)} = \bigg( \diffn{n}{}{x} + u_{n-1}\diffn{n-1}{}{x} + \cdots + u_1\diff{}{x} + u_0 \bigg) (y) \ . \]
Now, for any $g \in \matr{GL}_n\R$, the first row of $g \cdot \matr{W}(\gamma) := \matr{W}(\gamma)g$ determines another star-shaped curve $g \cdot \gamma = \gamma g$ whose components are a basis of solutions to $L_{\gamma}(y)=0$. This implies that:
\[ L_{g \cdot \gamma}(y) = L_{\gamma}(y) \ , \]
and hence, the vector $u := (u_0, \dots, u_{n-1})$ of coefficients of $L_{\gamma}(y)$ is uniquely determined by the orbit $[\gamma] := \matr{GL}_n\R \cdot \gamma$ of $\gamma$ induced by the free, transitive action of $\matr{GL}_n\R$ on $\matr{W}(\gamma)$, so we shall write $L_{[\gamma]}$, from now on. For any $\gamma \in [\gamma]$, then\footnote{For the equivalent approach via the calculus of jets, see e.g. \cite{Olver} and references therein.}, the parameter $x$ is its \textit{centro-affine arc length parameter}, $u$ is the (vector of) \textit{centro-affine curvatures}, and (\ref{FrenetSerret}) is the \emph{centro-affine Frenet-Serret equation} for the \emph{centro-affine Frenet-Serret frame} $\matr{W}(\gamma)$. Consequently, on the set of centro-affine curves over $I$:
\begin{equation}\label{PicardLindelof}
\mathcal{M} := \{ \gamma \in \mathcal{C}^{\infty}(I,\R^{1 \times n}) \ | \ \det\matr{W}(\gamma(x)) \neq 0 \, \forall \, x \in I\}/\matr{GL}_n\R \ ,
\end{equation}
it follows (by the above, and Picard-Lindel$\matr{\ddot{o}}$f's theorem) that the correspondence $\varphi:\mathcal{M} \rightarrow \mathcal{F}^n$ such that $\varphi[\gamma] = u$, is a bijection, and we may then endow $\mathcal{M}$ with the smooth structure of $\mathcal{F}^n$ via this correspondence.

To describe the tangent bundle of $\mathcal{M}$, we shall need an appropriate modification of the finite-dimensional notions of integral curve and vector field.
\begin{defn}
An \emph{admissible deformation $[\Gamma]$ of a centro-affine curve $[\gamma] := [\gamma](0)$} is the $\matr{GL}_n\R$-orbit of a smooth mapping $\Gamma \in \mathcal{C}^{\infty}(I' , \mathcal{M})$.
\end{defn}
In other words, $\Gamma = \{ \gamma(t) \}_{t \in I'}$ is a $t$-family of star-shaped curves $\gamma(t)$ in $\R^n$, which is smooth with respect to $x$ and $t$. Here, $I' \subset \R$ is a sufficiently small open interval around 0, which is determined by the extent to which $\gamma(t)$ is smooth in both variables. (Smoothness in $t$ amounts to the condition that $u(t)$ satisfy a certain PDE in $x$ and $t$; see (\ref{ScalarIntegrability}) and the comments following it.) Note that $t$ is not necessarily a centro-affine arc-length parameter for any given $x$.
\begin{defn}\label{EquivalentTangents}
Let $\Gamma_a = \{\gamma_a(t_a) \}_{t_a \in I_a}$ and $\Gamma_b = \{\gamma_b(t_b) \}_{t_b \in I_b}$ be smooth, 1-parameter families of star-shaped curves such that $[\Gamma_a]$ and $[\Gamma_b]$ are admissible deformations of the same centro-affine curve $[\gamma] := [\gamma_a(0)] = [\gamma_b(0)]$. \newline Let $u_a(x,t_a) := \varphi[\gamma_a(t_a)]$ and $u_b(x,t_b) := \varphi[\gamma_b(t_b)]$. Then $[\Gamma_a]$ and $[\Gamma_b]$ are \emph{equivalent at $[\gamma]$} whenever:
\[ \pdiff{u_a}{t_a}\bigg\vert_{t_a=0} = \pdiff{u_b}{t_b}\bigg\vert_{t_b=0} \ . \]
The tangent space $T_{[\gamma]}\mathcal{M}$ of $\mathcal{M}$ at $[\gamma]$ is the set of equivalence classes $\{ \pdiff{\left( \varphi[\gamma](t) \right)}{t}\big\vert_{t=0} \}$ of admissible deformations $[\gamma](t)$ of $[\gamma]$ which are equivalent at $[\gamma]$.
\end{defn}

Our main result is the following:
\begin{thm}\label{Theorem1}\hspace{0pt}
\begin{enumerate}
\item There is a smooth isomorphism $T_{[\gamma]}\mathcal{M} \cong \mathcal{F}^{1 \times n}$ of $\mathcal{F}$-modules. Modulo the kernel of the tangent vector mapping (see next), $T_{[\gamma]}\mathcal{M}$ can be identified with the cyclic left module $\mathcal{D}_x/\mathcal{D}_x L_{[\gamma]}$ over the ring $\mathcal{D}_x := \mathcal{F}[\del_x]$ of linear ordinary differential operators over $\mathcal{F}$.
\item Modulo the kernel of the tangent vector mapping, each $h \in \mathcal{F}^{1 \times n}$ is in correspondence with an $X \in \mathcal{D}_x$, uniquely determined by
\[ L_{[\gamma]} \big( X( [\gamma] ) \big) = -h \ . \]
Here, $X$ acts on $[\gamma]$ from the left by component-wise differentiation. The kernel of the tangent vector mapping is the space of solutions to the associated homogeneous equation $L_{[\gamma]} \big( X( [\gamma] ) \big) = 0$.
\item $[\gamma]$ allows an admissible deformation $[\gamma](t)$ \emph{iff} there is a 1-parameter family $[\gamma](t) \mapsto X(t)$, smooth in $t$, such that for each $t$, the evolution operator $\pdiff{}{t} - X(t)$ satisfies:
\begin{equation}\label{ScalarIntegrability}
L_{[\gamma](t)}\big( (\tfrac{\del}{\del t} - X(t))(y(t)) \big) = 0 
\end{equation}
for all solutions $y(t) \in \, <y_1(t), \dots, y_n(t)>_\R$.
\item Let $[\gamma_{\lambda}](t)$ be an admissible deformation of the centro-affine curve $[\gamma_{\lambda}]$, defined by the equation $L_{[\gamma](t)}(\psi(t)) = \lambda(t) \psi(t)$.\footnote{This is simply a 1-parameter family of linear ODEs, parameterized by $\lambda$, for each fixed $t$.} Then $\pdiff{\lambda}{t} \equiv 0$ \emph{iff} the curvature $u(t)$ satisfies an equation of the $n$-KdV hierarchy.
\end{enumerate}
\end{thm}
Theorems \ref{Theorem1}-(1) and \ref{Theorem1}-(2) are proven in $\S\S\,$\ref{Homogeneous}-\ref{Derivations}, using the homogeneous structure of $\mathcal{M}$.  An outline of the argument for Theorem \ref{Theorem1}-(3) is as follows: if $[\gamma](t)$ is an admissible deformation of $[\gamma]$, then for any $\gamma(t) \in [\gamma](t)$, $\matr{W}(\gamma(t))$ must satisfy $\tfrac{\del^2 \matr{W}(\gamma(t))}{\del x \del t} = \tfrac{\del^2 \matr{W}(\gamma(t))}{\del t \del x}$. Letting $\matr{X}(t) := \pdiff{\matr{W}(\gamma(t))}{t} \matr{W}(\gamma(t))^{-1}$, the equality of mixed derivatives of $\matr{W}(\gamma(t))$ is equivalent to the zero-curvature condition:
\begin{equation}\label{MatrixIntegrability}
\pdiff{\matr{F}}{t} - \pdiff{\matr{X}}{x} + [\matr{F}, \matr{X}] = 0 \ .
\end{equation}
It is straightforward to verify that if $\matr{X}(t)$ satisfies (\ref{MatrixIntegrability}) for each $t$, then $\matr{X}(t)$ is determined, up to solutions of the homogeneous version of (\ref{MatrixIntegrability}), by a $t$-family of differential operators for which (\ref{ScalarIntegrability}) holds. Note, as well, that for each $t$, $\matr{X}(t)$ is independent of the choice of $\gamma(t) \in [\gamma](t)$.

To prove that any evolution operator satisfying (\ref{ScalarIntegrability}) induces an admissible deformation of a centro-affine curve, it suffices to prove that over $\mathcal{F}_{x,t} := \mathcal{C}^{\infty}(I \times I', \R)$, the rank of $<\gamma(t), \gamma^{(1)}(t), \dots, \gamma^{(n-1)}(t) >_{\mathcal{F}_{x,t}}$ is equal to $n$, for all $t \in I'$ and all $\gamma(t) \in [\gamma](t)$. We explain this in $\S\,$\ref{AdmissibleDeformations}.

Regarding the iso-spectral condition of Theorem \ref{Theorem1}-(4), note that (\ref{ScalarIntegrability}) is equivalent to a PDE satisfied by $u(t)$, with the coefficients of $X$ as $\llq$parameters". Among the PDE satisfying (\ref{ScalarIntegrability}) are the equations $\pdiff{L}{t} = [P, L]$ of Lax type, for iso-spectral operators $P$ (described below). Geometrically, if $X_a$ and $X_b$ are operators corresponding to admissible deformations $\Gamma_a$ and $\Gamma_b$, then their induced flows commute \emph{iff} they are both iso-spectral deformations.

When $n=2$ and $3$, it is possible to explicitly compute the conditions on the (coefficients of the) operator $P$, and one finds that the curvature $u(t)$ must satisfy an equation of the $n$-KdV hierarchy (\cite{Pinkall, FujiKuro10, FujiKuro14, CIMB09, CIMB13, GLHMB}), but this becomes unfeasible for general $n$. Hence, to determine the analogous conditions for general $n$, we solve the same problem for the spectral-parameterized linear ODE $L(\psi)=\lambda\psi$ in $\S\,$\ref{IsospectralStrata}, using an approach reminiscent of \cite{Dickey, GelDik, Adler}.

Our approach can be summarized as follows: in the spectral-parameterized version of Theorem 1.3, a tangent vector $X_\lambda \in T_{[\gamma_\lambda]}\mathcal{M}$ is a formal Laurent series in $\lambda^{-1}$ with differential operator coefficients. The space of tangent vectors which are solutions to the homogeneous version of (\ref{MatrixIntegrability}), for $L(\psi)=\lambda\psi$, is essentially the space of power series in $\lambda^{-1}$, so all other tangent vectors must be polynomial in $\lambda$. (The latter such tangent vectors are generated from certain truncations of the former ones.) Further, such an $X_\lambda$ induces an iso-spectral flow -- i.e. $\pdiff{\lambda}{t} \equiv 0$ -- \emph{iff} the curvature $\varphi[\gamma_\lambda](t)$ satisfies an equation of the $n$-KdV hierarchy.

Thus, among the operators $P$ satisfying (\ref{ScalarIntegrability}), the ones which induce iso-spectral integrable deformations $[\gamma_\lambda](t)$ are precisely those of the form $P = X_\lambda\lvert_{\lambda=0}$, for $X_\lambda \in T_{[\gamma_\lambda]}\mathcal{M}$ which are polynomial-in-$\lambda$.


Some comments, on our assumptions and on the literature, are required. 

Regarding the usage of the centro-affine parameter $x$, let us recall (e.g. \cite{CIMB09, CQ, FujiKuro10, Pinkall, TerngWu15}) that if $\det\matr{W}(\gamma) > 0$ for all $x \in I$, then there is a unique re-parameterization of $I$ with smooth coordinate $s$ such that $\gamma(x) = \tilde{\gamma}(s)$ and $\det\matr{W}(\tilde{\gamma}(s)) = 1$ for all $s \in I$. We obtain a different operator $L_{\tilde{\gamma}}$, which necessarily has second-leading coefficient identically zero (by inspection of the co-factor expansion of $L_{\tilde{\gamma}}(y)$). The space of curves in question is then:
\[ \{ \tilde{\gamma} \in \mathcal{C}^{\infty}(I, \R^{1 \times n} )\ | \ \det\matr{W}(\tilde{\gamma}(s)) = 1 \, \forall \, s \in I \}/\matr{SL}_n\R \ , \]
and the orbit $\matr{SL}_n\R \cdot \tilde{\gamma}$ is an equivalence class of \emph{equi-centro-affine curves in $\R^n$}, with equi-centro-affine arc-length parameter $s$ and equi-centro-affine curvatures defined by the equi-centro-affine version of the Frenet-Serret equation for $\gamma$. For the purposes of this article, there is no particular advantage to working with equi-centro-affine curves over centro-affine curves, so we have elected to use the centro-affine parameter $x$ and centro-affine curvatures $u$.

It is worth noting that in the approaches mentioned above (except \cite{TerngWu15}), the r$\hat{\matr{o}}$le of the spectral parameter is largely absent. \cite{TerngWu15} is a continuation of earlier work on loop groups (e.g. \cite{TU00, TU11}), and makes use of results from \cite{DrinSok}, on the classification of KdV-type equations. However, we have found that many of their results can be proven directly from the method presented in this article, which differs from that of \cite{DrinSok}, and this was one motivation for this article.

As is well-known in integrable systems, the presence of the spectral parameter is key to extracting (some level of) explicit description of the evolution equations and their $\llq$hidden" symmetries. In our situation, the spectral parameter $\lambda$ allows us to introduce a 1-parameter family of equations such that the problem can be deduced to be solvable by general principles -- in this case, by Birkhoff factorization with respect to $\lambda$ of a fundamental solution matrix of the spectral-parameterized system. This is the link between the algebraic-analytical approach and the loop group geometric approach; a detailed analysis from the latter point-of-view is described in \cite{SegalWilson, PressSeg} and references therein (see also \cite{Palais, TU00}). The analysis of the particular consequences of this general observation, especially from the point-of-view of $\mathcal{D}$-modules, provided another stimulus for this article.

Let us finally mention that the analysis of the class of solutions of the Cauchy problem for $u(x,t)$ satisfying (\ref{MatrixIntegrability}), is a separate topic in itself, and detailed treatments can be found in the classical \cite{AbCl, BeDeTo}. (We will be content to accept that a solution can be constructed locally, according to Frobenius' theorem.)

\section{Admissible deformations of centro-affine curves}

\subsection{The homogeneous structure of $\mathcal{M}$}\label{Homogeneous}

To understand the homogeneous structure of $\mathcal{M}$, let us first observe the following:

\begin{prop}\label{GaugeTranslation}
For any star-shaped curves $\alpha$ and $\beta$, $\exists! \ P = \sum_{k=0}^{n-1} p_k\del_x^k \in \mathcal{D}_x$ such that $\beta = P(\alpha)$, where $P$ acts on $\alpha$ by component-wise differentiation.
\end{prop}

\prf The matrix $\matr{P} := \matr{W}(\beta)\matr{W}(\alpha)^{-1}$ satisfies the $1^{\matr{st}}$-order matrix linear ODE $\matr{P}^{(1)} - \matr{F}_{\beta}\matr{P} + \matr{P}\matr{F}_{\alpha} = 0$. By the special form of the matrices $\matr{F}_{\alpha}$ and $\matr{F}_{\beta}$, it is straightforward to check from this equation that the $k^\matr{th}$ row $P_{k-1}=(P_{k-1,0}, \dots, P_{k-1,n-1})$ of $\matr{P}$ is determined by the following conditions: $\matr{P}$ is invertible for all $x \in I$; if $k > 1$, then $P_{k-1}y_i = \big( P_0 y_i \big)^{(k-1)}$ for each $i$; and $P_0 y_i$ is a solution of $L_{[\beta]}(y)=0$ for each $i$. Hence, the desired differential operator is $P = \sum_{k=0}^{n-1} P_{0,k}\del_x^k$. \quad $\blacksquare$

\begin{defn}\label{Translations}
For fixed $[\alpha]$ and $[\beta]$,
\[ _{[\alpha]}\!V_{[\beta]} := \{ \matr{W}(\beta)g\matr{W}(\alpha)^{-1} \ | \ g \in \matr{GL}_n\R \} \]
is a groupoid, parameterizing all $\llq$translations" (gauge transformations) from $[\alpha]$ to $[\beta]$, corresponding to the translation of the respective curvatures $\varphi[\alpha] \mapsto \varphi[\beta]$. The group $_{[\gamma]}\!V_{[\gamma]}$ is the \emph{isotropy group of $[\gamma]$}.
\end{defn}

It will be helpful to understand the correspondence $\matr{W}(\beta)\matr{W}(\alpha)^{-1}$ independent of the representative of $[\alpha]$ and $[\beta]$. To that end, let us first describe $\mathcal{M}$ in the language of bundles. The Frenet-Serret frame $\matr{W}(\gamma) = \begin{bmatrix} y_1 & \cdots & y_n \end{bmatrix}$ is a globally smooth, non-zero section of the trivial principal $\matr{GL}_n\R$-bundle $\mathcal{P} := I \times \matr{GL}_n\R$. $\mathcal{P}$ is the associated frame bundle of the trivial vector bundle $E := I \times \R^{1 \times n}$, which carries the globally smooth, non-zero section $\gamma \in \Gamma(E)$. Since $\det\matr{W}(\gamma) \neq 0$ identically on $I$, $\Gamma(E) \cong < \gamma, \gamma^{(1)}, \dots, \gamma^{(n-1)} >_\mathcal{F}$, for any $\gamma \in [\gamma]$. Hence, $\mathcal{M}$ parameterizes all $\matr{GL}_n\R$-equivalence classes of bases of $\Gamma(E)$ of this form. 

It will also be useful to consider the affine space $\mathcal{A}(E)$ of connections on $E$. On $\mathcal{A}(E)$, there is a connection $_{\gamma}\!\nabla: \Gamma(E) \rightarrow \Omega^1 (E)$ for each star-shaped curve $\gamma$, whose contraction with $\diff{}{x}$ we denote by $_{\gamma}\!\nabla_x$. For notational convenience, we shall identify $\gamma^{(k)} \leftrightarrow \vect{e}_{k+1} \in \mathcal{F}^{n \times 1}$. Then for any $\vect{v} \in \Gamma(E) \cong <\gamma, \dots, \gamma^{(n-1)}>_\mathcal{F}$, $_{\gamma}\!\nabla_x$ acts on $\vect{v}$ from the left by:
\begin{equation}\label{Nabla_x}
_{\gamma}\nabla_x(\vect{v}) := ({\textstyle \diff{}{x}}\matr{I}_n + \matr{F^T})(\vect{v}) \ .
\end{equation}

Evidently, $_{g \cdot \gamma}\nabla = \, _{\gamma}\nabla$ for any $g \in \matr{GL}_n\R$, so we shall write $_{[\gamma]}\nabla$, from now on. Let us now consider $\Gamma(E)/\matr{GL}_n\R \cong < [\gamma], \dots, [\gamma^{(n-1)}]>_\mathcal{F}$.

\begin{prop}\label{Dmodules}
There is a left action of $\mathcal{D}_x$ on $\Gamma(E) \cong <\gamma, \dots, \gamma^{(n-1)}>_\mathcal{F}$, defined by $\del_x \cdot \sigma := \, _{[\gamma]}\!\nabla_x(\sigma)$ and extending $\mathcal{F}$-linearly for the rest of $\mathcal{D}_x$. Moreover, letting $L_{[\gamma]} = \del_x^n + \sum_{k=0}^{n-1} u_k \del_x^k$, the following isomorphism of left $\mathcal{D}_x$-modules holds for each $[\gamma]$:
\[ \mathcal{D}_{[\gamma]} := \mathcal{D}_x/\mathcal{D}_x L_{[\gamma]} \cong \, <[\gamma], \dots, [\gamma^{(n-1)}]>_\mathcal{F} \ . \]
\end{prop}

\prf The proof is by explaining the notation (see \cite{vdPS}). $\mathcal{D}_x = \mathcal{F}[\del_x]$ is the polynomial ring in the non-commutative indeterminate $\del_x$, which is defined to satisfy the Leibniz rule $\del_x \cdot f := f \del_x + f^{(1)}$ for all $f \in \mathcal{F}$ when considered as elements of $\mathcal{D}_x$. (Here, we still write $f^{(k)} := \diffn{k}{f}{x} \,$.) The ideal $\mathcal{D}_x L_{[\gamma]}$ is the principal left ideal generated by $L_{[\gamma]}$, so the quotient $\mathcal{D}_{[\gamma]} = \{ X + AL_{[\gamma]} \ | X, A \in \mathcal{D}_x \}$ is a left $\mathcal{D}_x$-module. An element of $\mathcal{D}_{[\gamma]}$ will be written $[X]_{[\gamma]}$, where we may always take $X \in [X]_{[\gamma]}$ to be the unique representative of degree $< n$ (by $\mathcal{D}_x$ a left Euclidean domain). $\mathcal{D}_{[\gamma]}$ is a cyclic left $\mathcal{D}_x$-module with generator $[1]_{[\gamma]}$, meaning that $\mathcal{D}_{[\gamma]} \cong \mathcal{D}_x \cdot [1]_{[\gamma]}$. It is also an $\mathcal{F}$-module of rank $n$, isomorphic to $\Gamma(E)$ under the identification of basis elements $[\del_x^k]_{[\gamma]} \leftrightarrow \gamma^{(k)}$. 
It is straightforward to check that the left action of $\del_x$ on $\mathcal{D}_{[\gamma]}$, and the left action of $_{[\gamma]}\!\nabla_x$, agree under this identification. 
It is also clear that $\gamma$ is a cyclic element of $\Gamma(E)$ with respect to $_{[\gamma]}\!\nabla_x$. Hence, we have an isomorphism of cyclic left $\mathcal{D}_x$-modules via the above left actions, and the identification of cyclic elements $[1]_{[\gamma]} \leftrightarrow [\gamma]$. \quad $\blacksquare$

\begin{prop}\label{Solutions}
The Frenet-Serret frame $\matr{W}(\gamma) = \begin{bmatrix} y_1 & \cdots & y_n \end{bmatrix}$ of a star-shaped curve $\gamma$ is the image of a linearly-independent (over $\R$) $n$-tuple $( \phi_1 , \dots , \phi_n )$ of left $\mathcal{D}_x$-linear homomorphisms $\phi_i:\mathcal{D}_{[\gamma]} \rightarrow \mathcal{F}^{n \times 1}$:
\[ \matr{W}(\gamma) = \begin{bmatrix} \phi_1([1]_{[\gamma]}) & \cdots & \phi_n([1]_{[\gamma]}) \end{bmatrix} \ . \]
Conversely if $\matr{\Phi} := \begin{bmatrix} \phi_1 & \cdots & \phi_n \end{bmatrix}$ is an $n$-tuple of left $\mathcal{D}_x$-linear homomorphisms $\phi_i:\mathcal{D}_{[\gamma]} \rightarrow \mathcal{F}^{n \times 1}$, then it is uniquely determined by its value at $[1]_{[\gamma]}$.
\end{prop}

\prf It suffices to explain the notion of left $\mathcal{D}_x$-linear maps (see \cite{vdPS}). A map $\phi: \mathcal{D}_{[\gamma]} \rightarrow \mathcal{F}$ is left $\mathcal{D}_x$-linear, or $\phi \in \matr{Hom}_{\mathcal{D}_x}(\mathcal{D}_{[\gamma]}, \mathcal{F})$ for short, if for all $X \in \mathcal{D}_x$ and for all $[A]_{[\gamma]} \in \mathcal{D}_{[\gamma]}$:
\[ X\big( \phi( [A]_{[\gamma]} ) \big) = \phi \big( X\cdot[A]_{[\gamma]} \big) \ . \]
Taking $A=1$ and letting $X$ be arbitrary, we see that any such $\phi$ is defined by its value at $[1]_{[\gamma]}$. 
Then taking $X=L_{[\gamma]}$, it follows that $y := \phi([1]_{[\gamma]})$ satisfies $L_{[\gamma]}(y)=0$. Moreover, $\phi([1]_{[\gamma]}) = 0$ \emph{iff} $\phi \equiv 0$, so the mapping $y \mapsto \phi_y$, taking a solution $y$ of $L_{[\gamma]}(y)=0$ to a left $\mathcal{D}_x$-linear map $\phi_y: \mathcal{D}_{[\gamma]} \rightarrow \mathcal{F}$ defined by $\phi_y([1]_{[\gamma]}) := y$, is also injective.

Hence, we have a vector space bijection between the space $<y_1, \dots, y_n>_\R$ of solutions to $L_{[\gamma]}(y)=0$, and the space $\matr{Hom}_{\mathcal{D}_x}(\mathcal{D}_{[\gamma]}, \mathcal{F}^{n \times 1})$ of all left $\mathcal{D}_x$-linear maps $\phi: \mathcal{D}_{[\gamma]} \rightarrow \mathcal{F}^{n \times 1}$. \quad $\blacksquare$

Having characterized $\mathcal{M}$ in terms of $\mathcal{D}_x$-modules, we can now give the following converse to Proposition \ref{GaugeTranslation}:

\begin{prop}\label{GaugeTranslationCor}
Let $[P]_{[\alpha]} := \sum_{k=0}^{n-1} p_k [\del_x^k]_{[\alpha]}$ such that $L_{[\beta]} \cdot [P]_{[\alpha]} = [0]_{[\alpha]}$ and
\[ \matr{P^T} := \begin{bmatrix} [P]_{[\alpha]} & \del_x \cdot [P]_{[\alpha]} & \cdots & \del_x^{n-1} \cdot [P]_{[\alpha]} \end{bmatrix} \in \mathcal{C}^{\infty}(I, \matr{GL}_n\R) \ , \]
where the columns of $\matr{P^T}$ are the components of the indicated coset (in order from the component of $[1]_{[\alpha]}$ in the top row, down to the component of $[\del_x^{n-1}]_{[\alpha]}$ in the bottom row). Then $\matr{P}\matr{W}(\alpha)$ is the Frenet-Serret frame of the star-shaped curve $P(\alpha) = (P(\alpha_1), \dots, P(\alpha_n)) \in [\beta]$.
\end{prop}

\prf By the form of $\matr{P}$ and the assumption $L_{[\beta]} \cdot [P]_{[\alpha]} = [0]_{[\alpha]}$, $\matr{P}$ satisfies $\matr{P}^{(1)} - \matr{F}_{\beta}\matr{P} + \matr{P}\matr{F}_{\alpha} = 0$. By Proposition \ref{Solutions}, $P(y_i) = \phi_i([P]_{[\alpha]})$, so $L_{[\beta]} \cdot [P]_{[\alpha]} = [0]_{[\alpha]}$ implies that $P(\alpha) \in [\beta]$, for any $\alpha \in [\alpha]$. \quad $\blacksquare$

For later use, let us observe the following:
\begin{cor}
$_{[\alpha]}V_{[\beta]}$ is isomorphic to the groupoid of all left $\mathcal{D}_x$-linear isomorphisms $\phi:\mathcal{D}_{[\beta]} \rightarrow \mathcal{D}_{[\alpha]}$. 
\end{cor}

\prf The left coset $[P]_{[\alpha]}$ of Proposition \ref{GaugeTranslationCor} has unique inverse $[P^+]_{[\beta]}$ in the sense that $[P^+]_{[\beta]} \cdot [P]_{[\alpha]} = [1]_{[\alpha]}$ and $[P]_{[\alpha]} \cdot [P^+]_{[\beta]} = [1]_{[\beta]}$. (These products are well-defined because $L_{[\beta]} \cdot [P]_{[\alpha]} = [0]_{[\alpha]}$ and $L_{[\alpha]} \cdot [P^+]_{[\beta]} = [0]_{[\beta]}$). Hence, the map $\phi_P: \mathcal{D}_{[\beta]} \rightarrow \mathcal{D}_{[\alpha]}$ defined by $\phi_P([A]_{[\beta]}) := [AP]_{[\alpha]}$ is a left $\mathcal{D}_x$-linear isomorphism.

Conversely, any left $\mathcal{D}_x$-linear isomorphism $\phi:\mathcal{D}_{[\beta]} \rightarrow \mathcal{D}_{[\alpha]}$, defining the coset $[P]_{[\alpha]} := \phi([1]_{[\beta]})$, necessarily satisfies $L_{[\beta]} \cdot [P]_{[\alpha]} = [0]_{[\alpha]}$. It can then be verified that its matrix representation $\matr{P}$, defined with respect to the cyclic bases of $\mathcal{D}_{[\beta]}$ and $\mathcal{D}_{[\alpha]}$, satisfies $\matr{P}^{(1)} - \matr{F}_{\beta}\matr{P} + \matr{P}\matr{F}_{\alpha} = 0$, so by Proposition \ref{GaugeTranslationCor}, it must be of the form $\matr{P} = \matr{W}(\beta)\matr{W}(\alpha)^{-1}$ for some $\alpha \in [\alpha]$ and $\beta \in [\beta]$. \quad $\blacksquare$

\subsection{Analytic first-order variations of the centro-affine curvature}\label{AnalyticVariations}

We now exam the relation between the space of direction vectors $\mathcal{F}^{1 \times n}$ and the tangent space $T_{[\gamma]}\mathcal{M}$ at $[\gamma]$, by determining analytic first-order variations in the centro-affine curvature $u \mapsto u + \epsilon h$ in terms of elements of $\mathcal{D}_{[\gamma]}$.

If $h = (h_0, \dots, h_{n-1}) \in \mathcal{F}^{1 \times n}$ is non-zero, then by Picard-Lindel$\matr{\ddot{o}}$f's theorem, the variation $u \mapsto u + \epsilon h$ induces an $\epsilon$-family of groupoids $V_{\epsilon}(h) := \, _{[\gamma]}V_{[\gamma_{\epsilon}(h)]}$ for all $\epsilon \in \R$, such that $\varphi[\gamma_{\epsilon}(h)] = u+\epsilon h$.

\begin{prop}\label{Variation}
For $h \neq 0$, every $\matr{P}_\epsilon(h) \in V_{\epsilon}(h)$ can be expanded in the form $\matr{P}_\epsilon(h)\matr{P}_0(h)^{-1} = \matr{I}_n + \epsilon \matr{X}(h) + o(\epsilon)$, where $\matr{P}_0(h) = \matr{P}_{\epsilon}(h)\vert_{\epsilon=0} \in \, _{[\gamma]}V_{[\gamma]}$, and $\matr{X}(h)$ is a solution of:
\begin{equation}\label{MatrixInhomogeneous}
\matr{X}^{(1)} - [ \matr{F}, \matr{X} ] = 
\begin{bmatrix}   0  & \cdots & 0       \\
              \vdots & \ddots & \vdots  \\
                  0  & \cdots & 0       \\
                -h_0 & \cdots & -h_{n-1} \end{bmatrix} \ .
\end{equation}
The matrix $\matr{X}(h)$ is necessarily of the same form as $\matr{P}$ described in Proposition \ref{GaugeTranslation}. The first row of $\matr{P}_{\epsilon}(h)\matr{P}_0(h)^{-1}$ induces an $\epsilon$-family $[P_\epsilon(h)]_{[\gamma]} = [1 + \epsilon X(h) + o(\epsilon)]_{[\gamma]}$, where $[X(h)]_{[\gamma]}$ is a solution of:
\begin{equation}\label{ScalarInhomogeneous}
L_{[\gamma]} \cdot [X(h)]_{[\gamma]} = - \sum_{k=0}^{n-1} h_k [\del_x^k]_{[\gamma]} \ .
\end{equation}
$\matr{X}(h)$ and $[X(h)]_{[\gamma]}$ are unique modulo solutions of their respective associated homogenous equations.
\end{prop}

\prf By Proposition \ref{GaugeTranslation}, $\matr{P}_{\epsilon}(h) =  \matr{W}(\gamma_{\epsilon}(h))\matr{W}(\gamma)^{-1}$ for some $\gamma_{\epsilon}(h) \in [\gamma_{\epsilon}(h)]$ and $\gamma \in [\gamma]$. Now, $\matr{P}_{\epsilon}(h)$ must satisfy:
\[ \matr{P}_{\epsilon}^{(1)} - [ \matr{F}, \matr{P}_{\epsilon}] + \epsilon (\vect{e}_n \otimes h)\matr{P}_{\epsilon} = 0 \ , \]
so by Picard-Lindel$\matr{\ddot{o}}$f's theorem with smooth parameters ($\S\,$2.5 of \cite{Teschl}), we may expand $\matr{P}_{\epsilon}(h) = \matr{P}_0(h) + \epsilon \matr{X}(h) + o(\epsilon)$, where $\matr{P}_0(h) \in \, _{[\gamma]}V_{[\gamma]}$ necessarily. Consequently, $\matr{X}(h)$ must satisfy for all $\epsilon$:
\[ \epsilon \big( \matr{X}^{(1)} - [\matr{F}, \matr{X}] + (\vect{e}_n \otimes h)\matr{P}_0 \big) = o(\epsilon) \ , \]
and since $\matr{P}_0(h)^{(1)} - [ \matr{F}, \matr{P}_0(h) ] = 0$, we see that $(\matr{X}\matr{P}_0^{-1})^{(1)} - [\matr{F}, \matr{X}\matr{P}_0^{-1}] = -\vect{e}_n \otimes h$, from which (\ref{MatrixInhomogeneous}) follows for $\matr{X}(h)\matr{P}_0(h)^{-1}$. 

By Proposition \ref{GaugeTranslation}, the first row of $\matr{P}_{\epsilon}(h)\matr{P}_0(h)^{-1}$ induces a left coset of the form $[P_{\epsilon}(h)]_{[\gamma]} = [1 + \epsilon X(h) + o(\epsilon)]_{[\gamma]}$, which must satisfy:
\begin{align*}
0 & = \left( L_{[\gamma]} + \epsilon \sum_{k=0}^{n-1} h_k \del_x^k \right) \cdot [ 1 + \epsilon X(h) + o(\epsilon) ]_{[\gamma]} \\
 & = \epsilon \left[ L_{[\gamma]} X(h) + \sum_{k=0}^{n-1} h_k \del_x^k \right]_{[\gamma]} + [o(\epsilon)]_{[\gamma]} \ ,
\end{align*}
for all $\epsilon$, from which (\ref{ScalarInhomogeneous}) follows. \quad $\blacksquare$

\begin{rmrk}\upshape
The space of solutions to the equation $\matr{X}^{(1)} - [\matr{F}, \matr{X}]=0$ will be denoted by
\[ \mathcal{E}[\gamma] = \{ \matr{W}(\gamma) m \matr{W}(\gamma)^{-1} \ | \ m \in \mathfrak{gl}_n\R \} \ . \]
(This goes by the name of $\llq$the eigenring of $L_{[\gamma]}$", in \cite{vdPS, Singer}.) The group of invertible elements of $\mathcal{E}[\gamma]$ is the isotropy group $_{[\gamma]}V_{[\gamma]}$.
\end{rmrk}

\begin{prop}\label{Injective}
The map $[X]_{[\gamma]} \mapsto L_{[\gamma]} \cdot [X]_{[\gamma]} = -[h]_{[\gamma]}$ is injective, for each $[X]_{[\gamma]} \in \mathcal{D}_{[\gamma]}/\mathcal{E}[\gamma]$.
\end{prop}

\prf For each $X \in \mathcal{D}_x$, there is a unique quotient-remainder pair $(\hat{X}, h)$ such that $L_{[\gamma]} X = \hat{X}L_{[\gamma]} - h$ and $h$ is of order $< n$, by $\mathcal{D}_x$ a left Euclidean domain. Now, if operators $X_1$ and $X_2$ are such that $X_1 - X_2 \in \mathcal{D}_xL_{[\gamma]}$, then they induce the same remainder $-h$, so to each $[X]_{[\gamma]}$, there is a unique $[h]_{[\gamma]}$ such that $L_{[\gamma]} \cdot [X]_{[\gamma]} = -[h]_{[\gamma]}$. Moreover, if $L_{[\gamma]} \cdot [X_1]_{[\gamma]} = L_{[\gamma]} \cdot [X_2]_{[\gamma]} = -[h]_{[\gamma]}$, then $[X_1 - X_2]_{[\gamma]} \in \mathcal{E}[\gamma]$, so $[X_1]_{[\gamma]} = [X_2]_{[\gamma]}$ as elements of $\mathcal{D}_{[\gamma]}/\mathcal{E}[\gamma]$. \quad $\blacksquare$

\begin{cor}\label{LieCorrespondence}
The following maps are well-defined and inverse to each other:
\[ \begin{array}{cccc}
V_1(h) & \longrightarrow & \mathcal{D}_{[\gamma]} / \mathcal{E}[\gamma] \\
\matr{P}_1(h) = \left( \matr{I}_n + \epsilon\matr{X}(h) + o(\epsilon)\right)\lvert_{\epsilon=1} & \mapsto & [X]_{[\gamma]} \ , 
\end{array} \]
\[ \begin{array}{cccc}
\mathcal{D}_{[\gamma]} / \mathcal{E}[\gamma] & \longrightarrow & V_1(h) \\
\ [X]_{[\gamma]}                               &    \mapsto      & \matr{P}_1(h) = \left( \matr{I}_n + \epsilon\matr{X}(h) + o(\epsilon)\right)\lvert_{\epsilon=1} \ , 
\end{array} \]
where $-[h]_{[\gamma]} = L_{[\gamma]} \cdot [X]_{[\gamma]}$ and $\matr{P}_{\epsilon}(h)= \matr{I} + \epsilon\matr{X}(h) + o(\epsilon) \in V_{\epsilon}(h)$ as in Proposition \ref{Variation}. In particular, $T_{[\gamma]}\mathcal{M} \cong \mathcal{D}_{[\gamma]}/\mathcal{E}[\gamma]$.
\end{cor}

\prf By Proposition \ref{Variation}, $\matr{X}(h) = \diff{}{\epsilon}\left( \matr{P}_{\epsilon}(h)\matr{P}_0(h)^{-1} \right)\lvert_{\epsilon = 0}$ is unique and well-defined modulo solutions to the homogeneous equation $\matr{X}^{(1)} = [\matr{F},\matr{X}]$. By the form of $\matr{X}(h)$, it is uniquely determined by its first row, whose components form the components of a coset $[X]_{[\gamma]}$ as in Proposition \ref{GaugeTranslationCor}.

By applying Proposition \ref{Variation} to Proposition \ref{Injective}, it follows that for each $[X]_{[\gamma]}$ such that $h \neq 0$, there is a unique $\matr{P}_{\epsilon}(h) \in V_{\epsilon}(h)$ such that $\matr{P}_{\epsilon}(h) = \matr{I}_n + \epsilon \matr{X}(h) + o(\epsilon)$, where $\matr{X}(h)$ is induced by $[X]_{[\gamma]}$ as in Proposition \ref{GaugeTranslationCor}. \ $\blacksquare$

\subsection{Derivations of $\mathcal{C}^{\infty}(\mathcal{M}, \R)$}\label{Derivations}

To show how the elements of $\mathcal{D}_{[\gamma]} / \mathcal{E}[\gamma]$ are tangent vectors in a more tangible way, we now consider the space of smooth functions $\mathcal{C}^{\infty}(\mathcal{M}, \R)$ and their infinitesimal variations.

Among the smooth functions from $\mathcal{M}$ to $\R$ are those which are pointwise induced by an element $\matr{P} \in \, _{[\gamma]}V_{[\gamma]}$ of the isotropy group of $[\gamma]$. Such functions exist because by Definition \ref{Translations}, $\matr{P} \in \, _{[\gamma]}V_{[\gamma]}$ \emph{iff} $\matr{P} = \matr{W}(\gamma)g\matr{W}(\gamma)^{-1}$ for some $g \in \matr{GL}_n\R$, so the characteristic polynomial $\det(s\matr{I}_n - \matr{P}) = s^n + \sum_{k=1}^{n} c_k s^{n-k}$ of $\matr{P}$ belongs to $\R[s]$. In particular, the $c_k = c_k(\matr{P})$ depend only on the centro-affine curve $[\gamma]$, and not on any particular star-shaped curve $\gamma \in [\gamma]$. Thus, a smooth function $\mathcal{S}: \mathcal{M} \rightarrow \R$ may be defined by $\mathcal{S}([\gamma]) := f(c_1, \dots, c_n)$ for some $f \in \mathcal{C}^{\infty}(\R^n, \R)$.

For $\mathcal{S} \in \mathcal{C}^{\infty}(\mathcal{M}, \R)$ not of the above type, we may use a similar construction with the characteristic polynomial of $\matr{P} = \matr{W}(\gamma) \matr{M} \matr{W}(\gamma)^{-1}$ for $\matr{M} \in \mathcal{C}^{\infty}(I, \matr{GL}_n\R)$, instead. This guarantees that the coefficients $c_k$ of $\det(x\matr{I} - \matr{M})$ are independent of the star-shaped representative $\gamma \in [\gamma]$, although now the coefficients depend on $x$, in general. Smooth functions $\mathcal{S}:\mathcal{M} \rightarrow \R$ may be defined as $\mathcal{S}([\gamma]) := \mathfrak{I}\big( f(c_1, \dots, c_n) \big)$ by using some $\R$-linear $\mathfrak{I}:\mathcal{F} \rightarrow \R$ and $f \in \mathcal{C}^{\infty}(\R^n, \R)$.

On the other hand, $\mathcal{S}:\mathcal{M} \rightarrow \R$ is smooth at each $[\gamma]$ \emph{iff} $\mathcal{S} \circ \varphi^{-1}: \mathcal{F}^n \rightarrow \R$ is. Thus, smooth functions $\mathcal{S}:\mathcal{M} \rightarrow \R$ may be defined as $\mathcal{S}([\gamma]) = \mathfrak{I}(F(\varphi[\gamma]))$, for some $F \in \mathcal{C}^{\infty}(\mathcal{F}^n,\mathcal{F})$. 

The G$\matr{\hat{a}}$teaux derivative $\dee\mathcal{S}_h$ of $\mathcal{S}$ along $h \in \mathcal{F}^{1 \times n}$ is defined for all $[\gamma] \in \mathcal{M}$ as:
\[ \dee\mathcal{S}_h[\gamma] := \lim_{\epsilon \rightarrow 0} \tfrac{1}{\epsilon} \big[ (\mathcal{S} \circ \phi^{-1} )(u + \epsilon h) - (\mathcal{S} \circ \phi^{-1} )(u) \big] \ . \]
We will assume that this generalizes the usual directional derivative, in the form:
\[ \dee\mathcal{S}_h[\gamma] = \mathfrak{I}\big( \hat{\nabla} \mathcal{S}_{[\gamma]}(h) \big) \ , \]
where the gradient $\hat{\nabla}\mathcal{S}_{[\gamma]} \in \mathcal{F}^{n \times 1}$ of $\mathcal{S}$ along $[\gamma]$ is defined by this, and from now on, we fix $\mathfrak{I}:\mathcal{F} \rightarrow \R$ to be an $\R$-linear mapping with kernel $\ker_\mathcal{F}\mathfrak{I}$ such that $\diff{}{x}\mathcal{F} \subseteq \ker_\mathcal{F}\mathfrak{I}$.

To extract a more tangible expression for $\hat{\nabla}\mathcal{S}_{[\gamma]}(h)$, let us identify $h = (h_0, \dots, h_{n-1})$ with $[h]_{[\gamma]} = \sum_{k=0}^{n-1} h_k [\del_x^k]_{[\gamma]}$, so that $[h]_{[\gamma]} = -L_{[\gamma]}\cdot [X(h)]_{[\gamma]}$ for $[X(h)]_{[\gamma]}$ as in Proposition \ref{Variation}. To ensure that $\dee\mathcal{S}_h[\gamma]$ is well-defined with respect to $[\gamma]$, we must then identify $\hat{\nabla} \mathcal{S}_{[\gamma]} := (\hat{\sigma}_0, \dots, \hat{\sigma}_{n-1})$ with $\hat{\sigma}_{[\gamma]} = \sum_{k=0}^{n-1} \hat{\sigma}_k \delta_k \in \mathcal{D}_{[\gamma]}^{\, *}$, where $\mathcal{D}_{[\gamma]}^{\, *}$ is the dual (over $\mathcal{F}$) of $\mathcal{D}_{[\gamma]}$:
\[ \mathcal{D}_{[\gamma]}^{\, *} = \matr{Hom}_\mathcal{F}\big( \mathcal{D}_{[\gamma]}, \mathcal{F} \big) \ . \]
This has the basis $<\delta_0, \dots, \delta_{n-1}>_\mathcal{F}$ dual to $<[1]_{[\gamma]}, \dots, [\del_x^{n-1}]_{[\gamma]}>_\mathcal{F}$:
\[ \delta_i\big( [\del_x^j]_{[\gamma]} \big) = \delta_{i,j} \quad \textrm{(Kronecker delta).} \]

We may then express $\hat{\nabla} \mathcal{S}_{[\gamma]}(h)$ by the following expression:
\begin{thm}\label{Gateaux}
The Gateaux derivative along $h$ is:
\[ \dee\mathcal{S}_h[\gamma] = -\mathfrak{I}\big[ \big( L_{[\gamma]}^* \cdot \hat{\sigma}_{[\gamma]} \big) \big( [X]_{[\gamma]} \big) \big] \ . \]
\end{thm}

The proof follows easily, after introducing the following notation and lemma:
\begin{defn}\label{LeftHorner}
For $A = \sum_{k=0}^{n} a_k \del_x^k$, define the finite sequence of operators $\{ H_k(A) \}_{k=0}^{n}$ by:
\[ \left\{ \begin{array}{ccl}
               H_0(A) & := & a_n \ , \\
               H_k(A) & := & H_{k-1}(A)\del_x + a_{n-k} \ . \end{array} \right. \]
\end{defn}
Since $A = H_k(A)\del_x^{n-k} + \sum_{j=0}^{n-k-1}a_j\del_x^j$ for all $k$ (by inspection), it follows that $H_k(A) = \sum_{j=0}^{k} a_{n-j} \del^{k-j}$ and $H_n(A) = A$.

Next, let us define $[\tilde{\gamma}] \in \mathcal{M}$ via (\ref{FrenetSerret}) by $L_{[\tilde{\gamma}]} := (-1)^n L_{[\gamma]}^*$, where $L_{[\gamma]}^*$ is the formal adjoint of $L_{[\gamma]}$:
\[ L_{[\gamma]}^* := \sum_{k=0}^{n} (-\del_x)^k u_k \ . \]
\begin{lemma}\label{LagrangeIdentity}
Associated to $[X]_{[\gamma]} \in \mathcal{D}_{[\gamma]}$ and $\sigma_{[\gamma]} \in \mathcal{D}_{[\gamma]}^{\, *}$, define the matrices $\matr{X}, \hat{\Sigma} \in \Gamma(\matr{End}(E^*))$ by
\[ \matr{X^T} = \begin{bmatrix} \, [X]_{[\gamma]} & \del_x \cdot [X]_{[\gamma]} & \cdots & \del_x^{n-1} \cdot [X]_{[\gamma]} \, \end{bmatrix} \ , \]
\[ \hat{\Sigma} = \begin{bmatrix} \, H_{n-1}(L_{[\gamma]})^* \cdot \hat{\sigma}_{[\gamma]} & \cdots & H_1(L_{[\gamma]})^* \cdot \hat{\sigma}_{[\gamma]} & \hat{\sigma}_{[\gamma]} \, \end{bmatrix} \ . \]
Then $\matr{tr}(\matr{X}\hat{\Sigma}) = \sum_{k=0}^{n-1} \big( H_{k}(L_{[\gamma]})^* \cdot \hat{\sigma}_{[\gamma]} \big)\big( \del_x^{n-1-k} \cdot [X]_{[\gamma]} \big)$, and:
\begin{equation}\label{LagrangeMagri}
\diff{}{x}\matr{tr}(\matr{X}\hat{\Sigma}) = \hat{\sigma}_{[\gamma]}\big( L_{[\gamma]} \cdot [X]_{[\gamma]} \big) - \big( L_{[\gamma]}^* \cdot \hat{\sigma}_{[\gamma]} \big)\big( [X]_{[\gamma]} \big) \ .
\end{equation}
\end{lemma}

\prf This can be verified directly, as in the proof of Proposition \ref{GaugeTranslation} (see also Corollary \ref{DualVariation}, below). Let us just comment on the conceptual reason for why: both $\matr{X}$ and $\hat{\Sigma}$ are representations of sections $\phi_X$ and $\phi_{\hat{\sigma}} \in \matr{End}_\mathcal{F}\mathcal{D}_{[\gamma]}^{\, *} \subset \Gamma(\matr{End}(E^*))$, respectively, so that as endomorphisms, the following relation holds:
\[ \diff{}{x}\matr{tr}(\phi_X \circ \phi_{\hat{\sigma}}) = \matr{tr}\big( (\del_x \cdot \phi_X) \circ \phi_{\hat{\sigma}} \big) + \matr{tr}\big( \phi_X \circ (\del_x \cdot \phi_{\hat{\sigma}}) \big) \ . \]
Since $\del_x$ acts on $\phi \in \matr{End}_\mathcal{F}\big( \mathcal{D}_{[\gamma]}^{\, *} \big)$ by $\del_x \cdot \phi := [ \, _{[\gamma]}\nabla_x^*, \phi ]$, the assertion follows by the form of $\matr{X}$ and $\hat{\Sigma}$. \quad $\blacksquare$

For more details, the reader is referred to the next section. The section after that returns to the consideration of admissible deformations.

\subsection{Dual Centro-Affine Curves in $\R^n$}\label{DualCurves}
In the following, we explore the dual versions of the results of $\S\S$\,\ref{Homogeneous}-\ref{AnalyticVariations}, and show the deeper meaning of the quotient $\hat{X}$ in the identity $L_{[\gamma]} X = \hat{X}L_{[\gamma]} - h$, in $\mathcal{D}$-module terms. (These observations will be applied in $\S\S$\,\ref{IsospectralStrata}-\ref{SpectralEigenring}.)

A fundamental solution matrix of the dual Frenet-Serret equation $\vect{y}^{(1)} = -\matr{F^T}\vect{y}$ is given by the dual Frenet-Serret frame  $\big(\matr{W}(\gamma)^{-1}\big)^\matr{T}$. Here, we express:
\[  \matr{W}(\gamma)^{-1} := \begin{bmatrix} \gamma^* & \gamma^{(1)*} & \cdots & \gamma^{(n-1)*} \end{bmatrix} = \begin{bmatrix} y_1^* \\ y_2^* \\ \vdots \\ y_n^* \end{bmatrix} \ , \]
where $\gamma^{(k)*} \in \Gamma(E^*) \cong \mathcal{F}^{n \times 1}$ and $y_k^* \in \mathcal{F}^{1 \times n}$, so that $\Gamma(E^*) \cong <\gamma^*, \dots, \gamma^{(n-1)*}>_\mathcal{F}$ is the $\mathcal{F}$-dual basis of $< \gamma, \dots, \gamma^{(n-1)}>_\mathcal{F}$, and $< y_1^*, \dots, y_n^*>_\R$ is the $\R$-dual basis of $<y_1, \dots, y_n>_\R$. (Note that $g \cdot \gamma = \gamma g$ implies that $g \cdot \gamma^* = g^{-1} \gamma^*$, for all $g \in \matr{GL}_n\R$.) Thus, $_{[\gamma]}\nabla_x^* = {\textstyle \diff{}{x}}\matr{I}_n - \matr{F}$ is the dual connection of $_{[\gamma]}\nabla_x$.

\begin{prop}\label{DualCyclic}
$\gamma^{(n-1) \, *}$ is the cyclic vector of $<\gamma^*, \dots, \gamma^{(n-1)*}>_\mathcal{F}$, so as left $\mathcal{D}_x$-modules:
\[ \matr{GL}_n\R \, \cdot <\gamma^*, \dots, \gamma^{(n-1)*}>_\mathcal{F} \ \cong \mathcal{D}_{[\gamma]}^{\, *} \cong \mathcal{D}_{[\tilde{\gamma}]} \ . \]
\end{prop}

\prf By the discussion before Proposition \ref{Dmodules}, and by definition of the dual basis $<\delta_0, \dots, \delta_{n-1}>_\mathcal{F}$, $\del_x$ acts on both $\Gamma(E^*) \cong <\gamma^*, \dots, \gamma^{(n-1)*}>_\mathcal{F}$ and $\mathcal{D}_{[\gamma]}^{\, *}$ from the left by $\del_x \cdot = \, _{[\gamma]}\nabla_x^*$. It is then straightforward to verify that:
\[ \delta_{n-1-k} = H_k(L_{[\gamma]})^* \cdot \delta_{n-1} \quad \textrm{and} \quad \gamma^{(n-1-k)*} = H_k(L_{[\gamma]})^* \big( \gamma^{(n-1)*} \big) \ , \]
and hence, that $L_{[\gamma]}^* \cdot \delta_{n-1} = 0$ and $L_{[\gamma]}^*\big( \gamma^{(n-1)*} \big) = 0$. Thus, $\delta_{n-1}$ and $\gamma^{(n-1)*}$ are cyclic vectors of their respective left $\mathcal{D}_x$-modules, so up to the choice of $\gamma \in [\gamma]$, the isomorphisms hold by identifying $\matr{GL}_n\R \cdot \gamma^{(n-1)*} \leftrightarrow \delta_{n-1} \leftrightarrow [1]_{[\tilde{\gamma}]}$.\quad $\blacksquare$

\begin{cor}\label{Phi}
There is a unique left $\mathcal{D}_x$-linear isomorphism $\Phi_{[\gamma]}: \mathcal{D}_{[\gamma]}^{\, *} \rightarrow \mathcal{D}_{[\tilde{\gamma}]}$ defined by $\Phi_{[\gamma]}(\delta_{n-1}) = [1]_{[\tilde{\gamma}]}$. This has standard matrix representation $\matr{S}_{[\gamma]}$ defined by $\big(\matr{W}(\gamma)^{-1}\big)^\matr{T} = \matr{S}_{[\gamma]}^\matr{T}\matr{W}(\tilde{\gamma})$, independent of the choice of $\gamma \in [\gamma]$.
The entries of $\matr{S}_{[\gamma]}$ are differential polynomials in the components of $u$, only. \ $\blacksquare$
\end{cor}

\begin{rmrk}\label{Lagrange}\upshape
There is an $\mathcal{F}$-bilinear perfect pairing $<~\cdot,~\cdot~>_{[\gamma]} : \mathcal{D}_{[\tilde{\gamma}]} \times \mathcal{D}_{[\gamma]} \rightarrow \mathcal{F}$, given by $< [A]_{[\tilde{\gamma}]} , [B]_{[\gamma]} >_{[\gamma]} := \Phi_{[\gamma]}^{-1}([A]_{[\tilde{\gamma}]})([B]_{[\gamma]})$. Geometrically, the left $\mathcal{D}_x$-linearity of $\Phi_{[\gamma]}$ implies that this pairing is preserved by $\del_x$, in the sense that
\[ \diff{}{x}< [A]_{[\tilde{\gamma}]} , [B]_{[\gamma]} >_{[\gamma]} = < \del_x \cdot [A]_{[\tilde{\gamma}]} , [B]_{[\gamma]} >_{[\gamma]} + < [A]_{[\tilde{\gamma}]} , \del_x \cdot [B]_{[\gamma]} >_{[\gamma]} \ . \]

Similarly, $\Phi_{[\gamma]}^*$ induces the perfect pairing $\mathcal{B}_{[\gamma]} : \mathcal{D}_{[\tilde{\gamma}]}^{\, *} \times \mathcal{D}_{[\gamma]}^{\, *} \rightarrow \mathcal{F}$ given by $\mathcal{B}_{[\gamma]}(Z,Y) := \Phi_{[\gamma]}^*(Z)(Y)$. It can be verified that for $Z$ and $Y$ of the form\footnote{Here, $<\tilde{\delta}_0, \dots, \tilde{\delta}_{n-1}>_\mathcal{F}$ is the basis of $\mathcal{D}_{[\tilde{\gamma}]}^{\, *}$ dual to $<[1]_{[\tilde{\gamma}]}, \dots, [\del_x^{n-1}]_{[\tilde{\gamma}]}>_\mathcal{F}$. } $Z = \sum_{k=0}^{n-1} z^{(k)} \tilde{\delta}_k$ and $Y = \sum_{k=0}^{n-1} y^{(k)} \delta_k$, where $y, z \in \mathcal{F}$ are arbitrary, $\mathcal{B}_{[\gamma]}(Z, Y)$ is equal to the bilinear concomitant of $L_{[\gamma]}$ applied to $z$ and $y$ \cite{Ince}:
\[ \mathcal{B}_{[\gamma]}(Z, Y) = \sum_{j=0}^{n-1} \sum_{k=0}^{j} y^{(j-k)} (-1)^k (u_{j+1}z)^{(k)} \ . \] 
For such $Z$ and $Y$, the property that $\mathcal{B}_{[\gamma]}$ is preserved by $\del_x$ reduces to the classical Lagrange identity (\emph{ibid.}):
\[ \diff{}{x}\mathcal{B}_{[\gamma]}(Z, Y) = z L_{[\gamma]}(y) - y L_{[\gamma]}^*(z) \ . \]
Analogous identities hold for $< [A]_{[\tilde{\gamma}]}, [B]_{[\gamma]} >_{[\gamma]}$ when $A = \sum_{k=0}^{n-1} H_{k}(L_{[\tilde{\gamma}]})^*(a) \del_x^{n-1-k}$ and $B = \sum_{k=0}^{n-1} H_{k}(L_{[\gamma]})^*(b) \del_x^{n-1-k}$ for arbitrary $a, b \in \mathcal{F}$.
\end{rmrk}

\begin{rmrk}\label{DualGaugeTranslation}\upshape
By Proposition \ref{GaugeTranslation} and Corollary \ref{Phi}, we see that if $\matr{P} = \matr{W}(\beta)\matr{W}(\gamma)^{-1}$, then $\matr{Q} := \big(\matr{P}^{-1}\big)^\matr{T}$ is a translation of the dual Frenet-Serret frame of $\gamma$
\[ \big(\matr{W}(\gamma)^{-1}\big)^\matr{T} = \begin{bmatrix} H_{n-1}(L_{[\gamma]})^*(\tilde{\gamma}) \\ \vdots \\ H_1(L_{[\gamma]})^*(\tilde{\gamma}) \\ \tilde{\gamma} \end{bmatrix} = \matr{S}_{[\gamma]}\matr{W}(\tilde{\gamma}) \] 
to the dual frame $\big(\matr{W}(\beta)^{-1}\big)^\matr{T}$ of $\beta$. Analogous to the form of $\matr{P}$ of Proposition \ref{GaugeTranslation}, the $k^\matr{th}$ row $Q_{k-1}=(Q_{k-1,0}, \dots, Q_{k-1,n-1})$ of $\matr{Q}$ is determined by: $\matr{Q}$ is invertible for all $x \in I$; if $k > 1$, then $Q_{n-1-k}\tilde{y}_i = H_k(L_{[\gamma]})^*\big( Q_{n-1} \tilde{y}_i \big)$ for each $i$; and $Q_{n-1}\tilde{y}_i$ is a solution of $\tilde{L}_{[\beta]}(y)=0$ for each $i$. Thus, $\matr{Q}$ is determined by its bottom row, which determines the element $Q_{[\gamma]} := \sum_{k=0}^{n-1} Q_{n-1,k}\delta_k$ of $\mathcal{D}_{[\gamma]}^{\, *}$, taking $\gamma^{(n-1)*}$ to $\beta^{(n-1)*} = Q_{[\gamma]}(\gamma^{(n-1)*})$, as well as the differential operator $\Phi_{[\gamma]}(Q_{[\gamma]})$ taking $\tilde{\gamma}$ to $\tilde{\beta} = \Phi_{[\gamma]}(Q_{[\gamma]})(\tilde{\gamma})$.
\end{rmrk}

\begin{cor}\label{DualVariation}
For nonzero $h = (h_0, \dots, h_{n-1}) \in \mathcal{F}^{1 \times n}$, if $\matr{P}_{\epsilon}(h) \in V_{\epsilon}(h)$, then $\matr{Q}_{\epsilon}(h) := \big( \matr{P}_{\epsilon}(h)^{-1} \big)^\matr{T}$ is of the form $\matr{Q}_{\epsilon}(h)\matr{Q}_0^{-1} = \matr{I}_n + \epsilon \matr{Y}(h) + o(\epsilon)$, where $\matr{Y} = -\matr{X^T}$ for $\matr{X}(h)$ as in Proposition \ref{Variation}. Expressing $-\matr{Y^T}=\matr{X}$ by its columns, we have the dual representation:
\[ -\matr{Y^T} = \begin{bmatrix} H_{n-1}(L_{[\gamma]})^* \cdot Y_{n-1} & \cdots & H_1(L_{[\gamma]})^* \cdot Y_{n-1} & Y_{n-1} \end{bmatrix} \ , \]
where the components of $H_k(L_{[\gamma]}^*)\cdot Y_{n-1}$ are ordered from the component of $\delta_0$ in the top row, down to the component of $\delta_{n-1}$ in the bottom row. Moreover, $\matr{Q}_{\epsilon}(h)\matr{Q}_0^{-1}$ induces the $\epsilon$-family of dual operators $Q_{\epsilon}(h) = 1 + \epsilon Y(h) + o(\epsilon) \in \mathcal{D}_{[\gamma]}^{\, *}$, where $Y(h)$ satisfies $L_{[\gamma]}^* \cdot Y(h) = h^* \cdot \delta_{n-1}$ and $h$ is identified with $\sum_{k=0}^{n-1} h_k \del_x^k$. $\matr{Y}(h)$ and $Y(h)$ are unique up to solutions of their respective associated homogeneous equations.
\end{cor}

\prf Expressing $\matr{X} = \begin{bmatrix} Y_0 & \cdots & Y_{n-1} \end{bmatrix}$ by its columns, the transpose of Equation (\ref{MatrixInhomogeneous}) is equivalent to the system of equations:
\[ \left\{ \begin{array}{ccl} -h_0 \vect{e}_n & = & \, _{[\gamma]}\nabla_x^* (Y_0) - u_0 Y_{n-1} \ , \\
      -h_k \vect{e}_n & = & \, _{[\gamma]}\nabla_x^* (Y_k) + Y_{k-1} - u_k Y_{n-1} \ , \quad 1 \leq k \leq n-1 \ .
   \end{array} \right. \]
Identifying $Y_k$ with the dual operator $Y_k = \sum_{j=0}^{n-1} Y_{j,k} \delta_j$, it follows by induction that for $1 \leq k \leq n-2$,
\[ Y_{n-1-k} = -\bigg( \sum_{j=0}^{k-1} (-\del_x)^{k-1-j} h_{n-1-j}\bigg)\cdot \delta_{n-1} + H_k(L_{[\gamma]})^* \cdot Y_{n-1} \ . \]
Hence, by $-h_0 \vect{e}_n = \, _{[\gamma]}\nabla_x^* (Y_0) - u_0 Y_{n-1}$, we find that $h^* \cdot \delta_{n-1} = L_{[\gamma]}^* \cdot Y_0$. \quad $\blacksquare$

\begin{cor}\label{ScalarOperatorDuality}
$X \in [X]_{[\gamma]}$ of Proposition \ref{Variation} and $Y \in \mathcal{D}_{[\gamma]}^{\, *}$ of Corollary \ref{DualVariation}, both induced by the same $\matr{P}_{\epsilon} \in V_{\epsilon}(h)$, are related by the operator identity:
\[ L_{[\gamma]}X = \Phi_{[\gamma]}(Y)^* L_{[\gamma]} - h \ , \]
where $\Phi_{[\gamma]}:\mathcal{D}_{[\gamma]}^{\, *} \rightarrow \mathcal{D}_{[\tilde{\gamma}]}$ is the map of Corollary \ref{Phi}.
\end{cor}

\prf For $X(h) \in [X(h)]_{\gamma}$ of minimal degree, $L_{[\gamma]} \cdot [X(h)]_{[\gamma]} = -[h]_{[\gamma]}$, as left cosets, \emph{iff} $L_{[\gamma]}X(h) = \hat{X}L_{[\gamma]} - h$, as operators, where $\hat{X}$ is uniquely determined by $L_{[\gamma]}$, $X(h)$, and the condition $\deg h < n$. Taking the formal adjoint, projecting $\mathcal{D}_x \twoheadrightarrow \mathcal{D}_{[\tilde{\gamma}]}$, and applying $\Phi_{[\gamma]}^{-1}$ shows that $L_{[\gamma]}^* \cdot \Phi_{[\gamma]}^{-1}(\hat{X}^*) = h^*\cdot\delta_{n-1}$. Hence, $Y(h)$ and $\Phi_{[\gamma]}^{-1}(\hat{X}^*)$ are both solutions of the inhomogeneous equation $L^* \cdot Y = h^* \cdot \delta_{n-1}$. By uniqueness in Proposition \ref{Variation} and Corollary \ref{DualVariation} with respect to $\matr{P}_{\epsilon}(h)$, they must coincide. \quad $\blacksquare$

\subsection{Admissible deformations and rank-preserving flows}\label{AdmissibleDeformations}

Let us introduce $\mathcal{F}_{x,t} := \mathcal{C}^{\infty}(I \times I', \R)$, where $I' \subseteq \R$ is equipped with the smooth parameter $t$, and let $\mathcal{D}_{x,t} := \mathcal{D}_x[\del_t]$ be the ring of linear partial differential operators in $\del_x$ and $\del_t$, which commute with each other: $\del_x \cdot \del_t = \del_t \cdot \del_x$. (We also write $f^{(k)}(x,t) = \pdiffn{k}{f(x,t)}{x}$ for $f \in \mathcal{F}_{x,t}$.) 

Let us recall the discussion following Theorem \ref{Theorem1}.3. Given a smooth 1-parameter family of star-shaped curves $\gamma(t)$, we wish to preserve the action of $\del_x$ on the space of sections $\Gamma(E)(t) \cong < \gamma(t), \gamma^{(1)}(t), \dots, \gamma^{(n-1)}(t)>_{\mathcal{F}_{x,t}}$ for each $t$, so we let $\del_x \cdot = \, _{[\gamma](t)}\nabla_x$ as in Equation (\ref{Nabla_x}). The action of $\del_t$ on $\vect{v} \in \Gamma(E)(t)$ is similarly given as $\del_t \cdot \vect{v} = \, _{[\gamma](t)}\nabla_t(\vect{v}) := \pdiff{\vect{v}}{t} + \matr{X}(x,t)\vect{v}$ for some $\matr{X} \in \mathcal{C}^{\infty}(I \times I', \mathfrak{gl}_n\R)$. However, 
\[ _{[\gamma](t)}\nabla_x \circ \, _{[\gamma](t)}\nabla_t = \, _{[\gamma](t)}\nabla_t \circ \, _{[\gamma](t)}\nabla_x \]
identically on $\Gamma(E)(t)$ \emph{iff} $\matr{F}$ and $\matr{X}$ satisfy Equation (\ref{MatrixIntegrability}).

This leads us to the present refinement of Theorem \ref{Theorem1}.2. Geometrically, this gives the necessary and sufficient conditions for integral curves in $\mathcal{M}$ to exist for some sufficiently small open interval $I'$, given in terms of smooth 1-parameter families of tangent vectors (as described by Theorem \ref{Theorem1}.1).

\begin{thm}\label{RankPreserving}
For all $t$, let $X(t) \in \mathcal{F}_{x,t}[\del_x]$. Then $[\gamma](t)$ is an admissible deformation of $[\gamma](0)$ \emph{iff} the equation $L_{[\gamma](t)} \cdot [\del_t - X(t)]_{[\gamma](t)} = [0]_{[\gamma](t)}$ holds over the left $\mathcal{D}_{x,t}$-module $\mathcal{D}_{[\gamma](t)} := \mathcal{D}_{x,t}/\mathcal{D}_{x,t} L_{[\gamma](t)}$.
\end{thm}

\prf Suppose that $[\gamma](t)$ is an admissible deformation of $[\gamma](0)$. Then by the discussion following Theorem \ref{Theorem1}, it followed that for all $y(t) \in \ker L_{[\gamma](t)}$, Equation (\ref{ScalarIntegrability}) is satisfied; that is,
\[ L_{[\gamma](t)}\big( \tfrac{\del y}{\del t}(t) - X(y)(t) \big) = 0 \ . \]
As in the proof of Proposition \ref{Solutions}, we may identify $y(t) = \phi([1]_{[\gamma](t)})$ for some $\phi \in \matr{Hom}_{\mathcal{D}_{x,t}}(\mathcal{D}_{[\gamma](t)}, \mathcal{F}_{x,t})$. Then by left $\mathcal{D}_{x,t}$-linearity of $\phi$,
\[ 0 = \phi\big( L_{[\gamma](t)} \cdot [\tfrac{\del}{\del t} - X(t)]_{[\gamma](t)} \big) \ . \]
Since $\phi$ is arbitrary, it follows that $L_{\gamma(t)} \cdot [\del_t - X(t)]_{[\gamma](t)} = [0]_{[\gamma](t)}$, as asserted.

To prove the converse, we must prove that if $L_{[\gamma](t)} \cdot [\del_t - X(t)]_{[\gamma](t)} = [0]_{[\gamma](t)}$, then the following left $\mathcal{D}_{x,t}$-module:
\[ \mathfrak{M} := \mathcal{D}_{x,t}/\mathcal{D}_{x,t}(L_{[\gamma](t)}, \del_t - X(t)) \subseteq \mathcal{D}_{[\gamma](t)} \]
is of rank $n$ over $\mathcal{F}_{x,t}$. Here, the left ideal in the quotient defining $\mathfrak{M}$ is:
\[ \mathcal{D}_{x,t}(L_{[\gamma](t)}, \del_t - X(t)) = \mathcal{D}_{x,t}L_{[\gamma](t)} + \mathcal{D}_{x,t}(\del_t - X(t)) \ . \]
Note that $\mathfrak{M}$ is a cyclic left $\mathcal{D}_{x,t}$-module \emph{iff} $[\del_t - X(t)]_{[\gamma](t)} = [0]_{[\gamma](t)}$ for all $t$.

To show that the rank condition is sufficient, let us assume that $\mathfrak{M}$ has rank $n$, for now. Then for any $\phi \in \matr{Hom}_{\mathcal{D}_{x,t}}(\mathfrak{M}, \mathcal{F}_{x,t}^{1 \times n})$, $\gamma(t) := \phi([1]_{[\gamma](t)})$ defines a 1-parameter family of invertible matrices $\matr{W}(\gamma(t))$ satisfying $\pdiff{\matr{W}(\gamma(t))}{x} = \matr{F}\matr{W}(\gamma(t))$, by Proposition \ref{Solutions}. 

We note that by Proposition \ref{Variation}, the zero-curvature condition (\ref{MatrixIntegrability}) for $\matr{F}$ is well-defined by the given $X(t)$: the matrix $\matr{X}$ induced by $X$ is unique up to elements of the eigenring $\mathcal{E}[\gamma](t)$, so any other matrix $\tilde{\matr{X}}$ in the tangent vector equivalence class of $\matr{X}$ is of the form $\tilde{\matr{X}}(t) = \matr{X}(t) + \matr{W}(\gamma(t))\matr{C}(t)\matr{W}(\gamma(t))^{-1}$, where $\matr{C} \in \mathcal{C}^{\infty}(I', \matr{GL}_n\R)$ (i.e. constant with respect to $x$). Then $\pdiff{(\tilde{\matr{X}} - \matr{X})}{x} = [ \matr{F}, \tilde{\matr{X}} - \matr{X}]$, and thus, the zero-curvature condition
\[ \pdiff{F(t)}{t} = \pdiff{\tilde{\matr{X}}(t)}{x} - [\matr{F}(t),\tilde{\matr{X}}(t)] = \pdiff{\matr{X}(t)}{x} - [\matr{F}(t), \matr{X}(t)] \]
is well-defined by the given $X(t)$, meaning that $[\gamma](t)$ is a well-defined admissible deformation of $[\gamma](0)$.

Let us now prove that $\mathfrak{M}$ is of rank $n$ over $\mathcal{F}_{x,t}$. For $A \in \mathcal{D}_{x,t}$, we shall write its coset in $\mathfrak{M}$ by $[A]_{\mathfrak{M}}$. Over $\mathfrak{M}$, we may use the relation $[\del_t]_\mathfrak{M} = [X]_\mathfrak{M}$ to recursively eliminate all instances of $\del_t$ in $[A]_\mathfrak{M}$, and then right-divide by $L_{[\gamma](t)}$ to obtain a unique expression of the form $[A]_\mathfrak{M} = \sum_{j=0}^{n-1} X_j[\del_x^j]_{\mathfrak{M}}$, for some $X_j \in \mathcal{F}_{x,t}$. This shows that $\mathfrak{M}$ is at most rank $n$ over $\mathcal{F}_{x,t}$, with spanning set $\{ [\del_x^j]_\mathfrak{M} \}_{j=0}^{n-1}$, so we next want to show that it is linearly-independent over $\mathcal{F}_{x,t}$.

We observe that for all $A \in \mathcal{D}_{x,t}$, there are unique $A' \in \mathcal{D}_{x,t}$ and $A'' \in \mathcal{D}_x$ such that $AL_{[\gamma](t)} = A'L_{[\gamma](t)}(\del_t-X) + A''L_{[\gamma](t)}$. 
This follows from the fact that if $X(t)$ such that $L_{[\gamma](t)} \cdot [\del_t - X(t)]_{[\gamma](t)} = [0]_{[\gamma](t)}$, then there is a unique $\hat{X}$ such that $L_{[\gamma](t)}(\del_t - X) = (\del_t - \hat{X})L_{[\gamma](t)}$. 
Hence, we obtain the asserted re-expression of the operator $AL_{[\gamma](t)}$, by moving one power of $\del_t$ in $A$ from the left side of $L_{[\gamma](t)}$ to its right side according to $\del_t L_{[\gamma]}(t) = L_{[\gamma](t)}(\del_t - X) + \hat{X}L_{[\gamma(t)]}$.

Finally, we show that if there are $a_0, \dots, a_{n-1} \in \mathcal{F}_{x,t}$ and $A, B \in \mathcal{D}_{x,t}$ such that $a_0 + \cdots + a_{n-1} \del_x^{n-1} = AL_{[\gamma](t)} + B(\del_t - X)$, then all of the $a_j$ are zero. By the above observation regarding $AL_{[\gamma](t)}$, there are unique $A' \in \mathcal{D}_{x,t}$ and $A'' \in \mathcal{D}_x$ such that:
\[ \sum_{j=0}^{n-1} a_j \del_x^j - A''L_{[\gamma](t)} = (A'L_{[\gamma](t)} + B)(\del_t - X) \ . \]
The left-hand side belongs to $\mathcal{D}_x$, which forces $B = -A'L_{[\gamma](t)}$. By comparing degrees of $\del_x$ in $\sum_{j=0}^{n-1} a_j \del_x^j = A''L_{[\gamma](t)}$, we see that $A''=0$ and $a_j = 0$. \quad $\blacksquare$

\subsection{Admissible Deformations according to Isospectrality}\label{IsospectralStrata}

Let $[\gamma] \in \mathcal{M}$ have curvature $\varphi[\gamma] = u$, and let $[\gamma_{\lambda}] \in \mathcal{M}$ be a 1-parameter family of centro-affine curves defined via the Frenet-Serret equation (\ref{FrenetSerret}) for the ODE $L_{[\gamma]}(\psi) = \lambda\psi$, so that $\varphi[\gamma_{\lambda}] = u - (\lambda, 0, \dots, 0)$. Here, the $\llq$spectral parameter" $\lambda$ is independent of the centro-affine parameter $x \in I$, and we may re-write (\ref{FrenetSerret}) in this case as
\begin{equation}\label{SpectralFrenetSerret}
\pdiff{}{x}\matr{W}(\gamma_{\lambda}) = (\matr{F} + \lambda\matr{E}_{n,1})\matr{W}(\gamma_{\lambda}) \ ,
\end{equation}
where $\matr{F}$ is the matrix for the Frenet-Serret equation of $[\gamma]$, and $\matr{E}_{n,1}$ is the matrix with $1$ in the $(n,1)$ entry, and $0$ elsewhere.

Our goal, here, is to prove Theorem \ref{Theorem1}.4: for an admissible deformation $[\gamma_{\lambda}](t)$ of $[\gamma_{\lambda}]$, the curvature vector $u(t)$ of $[\gamma](t) = [\gamma_{\lambda=0}](t)$ satisfies an equation of the $n$-KdV hierarchy \emph{iff} $\pdiff{\lambda}{t} \equiv 0$ for all $t$. To do so, we first state the spectral parameterized version of Theorem \ref{RankPreserving} over the extended ring $\mathcal{D}_{x,t}^\lambda := \mathcal{D}_{x,t}((\lambda^{-1}))$ (formal Laurent series in $\lambda^{-1}$), where $\lambda$ only commutes with $\del_x$, \emph{a priori}. Let us denote 
\[ \mathcal{D}_{[\gamma_{\lambda}](t)} := \mathcal{D}_{x,t}^\lambda / \mathcal{D}_{x,t}^\lambda L_{[\gamma_{\lambda}](t)} \]
and 
\[ \mathcal{E}[\gamma_{\lambda}](t) := \{ [A]_{[\gamma_{\lambda}](t)} \in \mathcal{D}_{[\gamma_{\lambda}](t)} \ | \ L_{[\gamma_{\lambda}](t)} \cdot [A]_{[\gamma_{\lambda}](t)} = [0]_{[\gamma_{\lambda}](t)} \} \ . \]

\begin{prop}\label{SpectralizedRankPreserving}
$[\gamma_{\lambda}](t)$ is an admissible deformation of $[\gamma_{\lambda}](0)$ \emph{iff} for all $t$, $[\del_t - X]_{[\gamma_{\lambda}](t)} \in \mathcal{E}[\gamma_{\lambda}](t)$. (Here, $X \in \mathcal{F}_{x,t}[\del_x]((\lambda^{-1}))$.) \ $\blacksquare$
\end{prop}

We now determine which admissible deformations correspond to isospectral deformations:
\begin{thm}\label{Isospectral}
Let $[\del_t-P]_{[\gamma_{\lambda}](t)} \in \mathcal{E}[\gamma_{\lambda}](t)$, where $[P]_{[\gamma_{\lambda}](t)} \in \mathcal{F}_{x,t}[\del_x][\lambda]$. 
Then $\pdiff{\lambda}{t} \equiv 0$ for all $t$ \emph{iff} $[P]_{[\gamma_{\lambda}](t)} \in \mathcal{D}_{[\gamma_{\lambda}](t)} / \mathcal{E}[\gamma_{\lambda}](t)$.
\end{thm}

\prf The condition $L_{[\gamma_{\lambda}](t)} \cdot [\del_t - P]_{[\gamma_{\lambda}](t)} = [0]_{[\gamma_{\lambda}](t)}$ implies that there is a unique $Q \in \mathcal{F}_{x,t}[\del_x][\lambda]$ such that $L_{[\gamma_{\lambda}](t)}(\del_t - P) = (\del_t - Q)L_{[\gamma_{\lambda}](t)}$. As a result, we may re-write this over $\mathcal{D}_{[\gamma](t)} = \mathcal{D}_{x,t}^\lambda/\mathcal{D}_{x,t}^\lambda L_{[\gamma](t)}$ as:
\[ \left[ \pdiff{(L_{[\gamma](t)} - \lambda)}{t} \right]_{[\gamma](t)} + L_{[\gamma](t)} \cdot [P]_{[\gamma](t)} = \lambda[P-Q]_{[\gamma](t)} \ . \]
Since $P = \sum_{j=0}^{N} P_j \lambda^{N-j}$, it follows that $Q = \sum_{j=0}^{N} Q_j \lambda^{N-j}$. Therefore, $\pdiff{\lambda}{t} = 0$ \emph{iff} the following system of equations holds:
\[ \left\{ \begin{array}{cccl}
0 & = & [P_0 - Q_0]_{[\gamma](t)} \ , & \\
L_{[\gamma](t)} \cdot [P_j]_{[\gamma](t)} & = & [P_{j+1} - Q_{j+1}]_{[\gamma](t)} \ , & 0 \leq j \leq N-1 \ , \\
\left[ \pdiff{ L_{[\gamma](t)} }{t} \right]_{[\gamma](t)} + L_{[\gamma](t)} \cdot [P_N]_{[\gamma](t)} & = & [0]_{[\gamma](t)} \ . & \end{array} \right. \]
The first $N$ equations imply that $[P]_{[\gamma_{\lambda}](t)}$ is a linear combination of elements of $\mathcal{E}[\gamma_{\lambda}](t)$ (see $\S\,$\ref{SpectralEigenring} for details), then $\llq$truncated" in the quotient space $\mathcal{D}_{[\gamma_{\lambda}](t)} / \mathcal{E}[\gamma_{\lambda}](t)$ of non-trivial tangent vectors. (The meaning of the final equation is described, next.) \quad $\blacksquare$

\begin{cor}
For $[\del_t-P]_{[\gamma_{\lambda}](t)} \in \mathcal{E}[\gamma_{\lambda}](t)$ with $[P]_{[\gamma_{\lambda}](t)} \in \mathcal{F}_{x,t}[\del_x][\lambda]$, the curvature vector $u(t)$ of $[\gamma](t) = [\gamma_{\lambda=0}](t)$ satisfies an equation of the $n$-KdV hierarchy.
\end{cor}

\prf $\pdiff{\lambda}{t} \equiv 0$ \emph{iff} the curvature vector $u(t)$ satisfies an equation of the type described in the proof of Theorem \ref{Isospectral}. The $n$-KdV hierarchy is conventionally defined using elements of the kernel of the Adler map $A^{(\lambda)}$ for $L_{[\gamma_{\lambda}]}$ (e.g. \cite{Dickey}), but as shown next in Proposition \ref{GeneratorIsomorphism}, $\mathcal{E}[\gamma_{\lambda}]$ and the kernel of $A^{(\lambda)}$ are isomorphic as $\R((\lambda^{-1}))$-algebras, so the construction of the $n$-KdV hierarchy of $L_{[\gamma]} - \lambda$ may also be performed in the eigenring $\mathcal{E}[\gamma_{\lambda}]$, and the assertion follows. \quad $\blacksquare$

\begin{prop}\label{GeneratorIsomorphism}
Let $A^{(0)}(X)$ denote the \emph{Adler mapping} \cite{Dickey, Adler} of $X$ with respect to $L_{[\gamma]} = \del_x^n + \sum_{k=0}^{n-1} u_k \del_x^k$:
\[ A^{(0)}(X) := L_{[\gamma]}(XL_{[\gamma]})_+ - (L_{[\gamma]}X)_+ L_{[\gamma]} \ . \]
Then each pseudodifferential operator $X = \del^{-1}X_0 + \cdots + \del^{-n}X_{n-1}$ annihilated by $A^{(0)}$ uniquely corresponds to the element $Q \in \mathcal{D}_{[\gamma]}^{\, *}$ defined by:
\[ Q := \sum_{j=0}^{n-1} X_j \delta_j = \sum_{j=0}^{n-1} (X_{n-1-j} H_j(L_{[\gamma]})^*)\cdot \delta_{n-1} \ , \]
satisfying $L_{[\gamma]}^* \cdot Q = 0$.
\end{prop}

\prf We first observe that for $X$ of the given form:
\begin{align*}
(L_{[\gamma]}X)_+ & = \sum_{j=1}^{n} \sum_{k=j}^{n} u_k \del^{k-j} X_{j-1} = \sum_{j=0}^{n-1} \sum_{k=n-j}^{n} u_k \del^{k-n+j} X_{n-1-j} \\
 & = \sum_{j=0}^{n-1} \sum_{k=0}^{j} u_{n-k} \del^{j-k} X_{n-1-j} = \sum_{j=0}^{n-1} H_{n-1-j}(L_{[\gamma]}) X_j \ .
\end{align*}
Hence, for $Q$ as in the statement of the proposition, it follows from Corollary \ref{Phi} that $[(L_{[\gamma]}X)_+^*]_{[\tilde{\gamma}]} = [\Phi_{[\gamma]}(Q)]_{[\tilde{\gamma}]}$. Consequently:
\[ 0 = A^{(0)}(X) \quad \Rightarrow \quad L_{[\gamma]}^* \cdot [\Phi_{[\gamma]}(Q)]_{[\tilde{\gamma}]} = [0]_{[\tilde{\gamma}]} \ , \]
which implies that $Q \in \mathcal{D}_{[\gamma]}^{\, *}$ such that $L_{[\gamma]}^* \cdot Q = 0$, by left $\mathcal{D}$-linearity of $\Phi_{[\gamma]}$. 

Conversely, the above computation shows that any $Q = \sum_{j=0}^{n-1} X_j \delta_j$ uniquely defines the operator $\Phi_{[\gamma]}(Q)$ of degree $< n$ from the left coset of $[\Phi_{[\gamma]}(Q)]_{[\tilde{\gamma}]}$, such that $(L_{[\gamma]}X)_+^* = \Phi_{[\gamma]}(Q)$ for the pseudodifferential operator $X := \del^{-1} \sum_{j=0}^{n-1} \del^{-j} X_j$. To show that $A^{(0)}(X)=0$ whenever $Q$ is left-annihilated by $L_{[\gamma]}^*$, we note that for such a $Q$, there is a unique differential operator $P$ of degree $< n$ such that $L_{[\gamma]}P = (L_{[\gamma]}X)_+L_{[\gamma]}$. Hence, we wish to show that $P = (XL_{[\gamma]})_+$.

By the same reasoning as in the above, but applied to $L_{[\gamma]}P = (L_{[\gamma]}X)_+ L_{[\gamma]}$, it follows that $P = \Phi_{[\tilde{\gamma}]}(\tilde{P})$ for some $\tilde{P} \in \mathcal{D}_{[\tilde{\gamma}]}^{\, *}$ which defines a unique pseudodifferential operator $Y$ of the form $Y = \del^{-1}Y_0 + \cdots + \del^{-n}Y_{n-1}$, and $P = (YL_{[\gamma]})_+$. Using the differential-integral decomposition of pseudo-differential operators, we see that:
\begin{align*}
L_{[\gamma]}(YL_{[\gamma]})_+ = (L_{[\gamma]}X)_+ L_{[\gamma]} & = L_{[\gamma]}XL_{[\gamma]} - (L_{[\gamma]}X)_- L_{[\gamma]} \\
 & = L_{[\gamma]}(XL_{[\gamma]})_+ + L_{[\gamma]}(XL_{[\gamma]})_- - (L_{[\gamma]}X)_- L_{[\gamma]}
\end{align*}
\[ \Leftrightarrow \quad L_{[\gamma]}\big( (Y-X)L_{[\gamma]} \big)_+ = L_{[\gamma]}(XL_{[\gamma]})_- - (L_{[\gamma]}X)_- L_{[\gamma]} \ . \]
The left-hand side is a differential operator of degree $\leq 2n-1$, and the right-hand side is the difference of two pseudodifferential operators of degree $\leq n-1$. For this equality to hold, it is necessary that the top $n$ terms of the left-hand side vanish identically. By expanding $\big((Y-X)L_{[\gamma]}\big)_+$ and inspecting the triangular system of equations satisfied by $Y_k - X_k$ for each $k$, it can be verified that $Y_k = X_k$, and thus, $P = (XL_{[\gamma]})_+$, as desired. As a result, we conclude that if $L_{[\gamma]}^* \cdot Q = 0$, then $A^{(0)}(X)=0$. \quad $\blacksquare$

\begin{cor}
The kernel of the Adler mapping $A^{(\lambda)}$ of $L_{[\gamma_{\lambda}]}$ is isomorphic, as an $\R((\lambda^{-1}))$-algebra, to $\mathcal{E}[\gamma_{\lambda}]$.
\end{cor}

\prf By Corollary \ref{ScalarOperatorDuality}, we have an isomorphism between elements $Q \in \mathcal{D}_L^{\, *}$ satisfying $L^* \cdot Q = 0$, and elements $[P]_L \in \mathcal{D}_L$ satisfying $L \cdot [P]_L = [0]_L$. Thus, by Proposition \ref{GeneratorIsomorphism}, we have an isomorphism of $\R$-algebras between the kernel of $A^{(0)}$ and $\mathcal{E}[\gamma]$. Evidently, Proposition \ref{GeneratorIsomorphism} holds for $L_{[\gamma_{\lambda}]} = L_{[\gamma]} - \lambda$, after extending $\mathcal{D}_x$ to $\mathcal{D}_x^\lambda$ and $\R$ to $\R((\lambda^{-1}))$. \quad $\blacksquare$

We finish this section by showing that for different iso-spectral admissible deformations (labelled by different $\llq$time" variables) of the same centro-affine curve, the flows commute with each other.

By Theorem \ref{CanonicalBasis}, we shall use the following basis for $\mathcal{E}[\gamma_{\lambda}]$:
\[ \mathcal{E}[\gamma_{\lambda}] = < [1]_{[\gamma_{\lambda}]}, [\Omega]_{[\gamma_{\lambda}]}, \dots, [\Omega^{n-1}]_{[\gamma_{\lambda}]} >_{\R((\lambda^{-1}))} \ , \]
where $\Omega = \sum_{k=0}^{\infty} \Omega_k \lambda^{-k}$ is uniquely determined by the conditions: $\deg \Omega_k < n-1$; $\Omega_0 = \del_x + \tfrac{1}{n}u_{n-1}$; and $[\Omega^n]_{[\gamma_{\lambda}]} = [\lambda]_{[\gamma_{\lambda}]}$. (For more details, see the following section.) Let $\bar{\Omega}$ be the matrix representation of $[\Omega]_{[\gamma_{\lambda}]}$, as in Proposition \ref{GaugeTranslationCor}.

\begin{lemma}\label{ConservationLaws}
Let $[\del_t - P]_{[\gamma_{\lambda}]} \in \mathcal{E}[\gamma_{\lambda}](t)$ such that $\pdiff{\lambda}{t} \equiv 0$, and let $\matr{P}$ be the matrix representation of $[P]_{[\gamma_{\lambda}](t)}$, as in Proposition \ref{GaugeTranslationCor}, so that $_{[\gamma_{\lambda}](t)}\nabla_t \cdot = \del_t + \matr{P}$. Then $[ \ _{[\gamma_{\lambda}](t)}\nabla_t, \bar{\Omega}^k] = 0$ for all $k$, and consequently:
\[ [ \tfrac{\del}{\del t}\left( \Omega^k \right) ]_{[\gamma_{\lambda}](t)} = [ \, [P, \Omega^k] \, ]_{[\gamma_{\lambda}](t)} \ . \]
\end{lemma}

\prf By Theorem \ref{RankPreserving}, $_{[\gamma_{\lambda}](t)}\nabla_t$ commutes with $_{[\gamma_{\lambda}](t)}\nabla_x$. Since $\mathcal{E}[\gamma_{\lambda}](t)$ is an associative algebra, $[ \ _{[\gamma_{\lambda}](t)}\nabla_t, \bar{\Omega}^k]$ must also commute with $_{[\gamma_{\lambda}](t)}\nabla_x$, for each $k$. But the first column of $[ \ _{[\gamma_{\lambda}](t)}\nabla_t, \bar{\Omega}^k ]$ corresponds to $[ \, [\del_t - P, \Omega^k] \, ]_{[\gamma_{\lambda}](t)}$, so this is an element of $\mathcal{E}[\gamma_{\lambda}](t)$. We claim that $[ \, [\del_t - P, \Omega^k] \, ]_{[\gamma_{\lambda}](t)} = [0]_{[\gamma_{\lambda}](t)}$ for all $k$.

To show this, first recall from the proof of Theorem \ref{Isospectral} that
\[ \left[ \pdiff{ L_{[\gamma](t)} }{t} \right]_{[\gamma](t)} + L_{[\gamma](t)} \cdot [P_N]_{[\gamma](t)} = [0]_{[\gamma](t)} \ . \]
Now, $[P]_{[\gamma_{\lambda}](t)}$ is a polynomial-in-$\lambda$ truncation of $\lambda^N$ times a linear combination over $\R((\lambda^{-1}))$ of the $[\Omega^k]_{[\gamma_{\lambda}]}$. Using Equation (\ref{SequenceX}) (see $\S\,$\ref{SpectralEigenring}), we may use $\llq$un-truncated" linear combinations of the $[\Omega^k]_{[\gamma_{\lambda}]}$ to define the operator $P_{N+1} - Q_{N+1}$ by (cf. proof of Theorem \ref{Isospectral})
\[ L_{[\gamma](t)} \cdot [P_N]_{[\gamma](t)} =: [P_{N+1} - Q_{N+1}]_{[\gamma](t)} \ , \]
so that
\[ \left[ \pdiff{ L_{[\gamma](t)} }{t} \right]_{[\gamma](t)} = - [P_{N+1} - Q_{N+1}]_{[\gamma](t)} \ . \]
The conclusions of $\S\,$\ref{SpectralEigenring} -- in particular, Remarks \ref{SecondRemark} and Equation (\ref{y_j}) -- guarantee that this cannot have any constants (either in $x$ or $t$) in the leading $\lambda^0$ term, so $[ \, [\del_t, \Omega^k ] \, ]_{[\gamma_{\lambda}](t)}$ contains no constants for $k \neq n\Z$. To show that $[ \, [P, \Omega^k] \, ]_{[\gamma](t)}$ has no constants, we observe that the only constant terms of $[P]_{[\gamma_{\lambda}](t)}$ and $[\Omega^k]_{[\gamma_{\lambda}](t)}$ are the leading terms, and these vanish in the commutator, leaving no constant terms elsewhere. Therefore, by Lemma \ref{Step1}, $[ \, [\del_t - P, \Omega^k] \, ]_{[\gamma_{\lambda}](t)} = [0]_{[\gamma_{\lambda}](t)}$. \quad $\blacksquare$ \\

Thus far, we have fixed the deformation variable $t$, but it is clear that we may simultaneously consider multiple $\llq$time" variables. As such, let $t_i$ and $t_j$ be any two deformation variables describing admissible deformations of $[\gamma_{\lambda}]$, such that $\del_{t_i}$ and $\del_{t_j}$ all commute with $\del_x$ and with each other. Then the next proposition establishes that as long as these deformations are iso-spectral, the flows commute with each other (on the space of sections over $[\gamma_{\lambda}]$). Let 
\[ \mathcal{D}_{[\gamma_{\lambda}](t_i,t_j)} = \mathcal{D}_{x,t_i,t_j}/\mathcal{D}_{x,t_i,t_j}L_{[\gamma_{\lambda}](t_i,t_j)} \ . \]
\begin{prop}\label{CommutingFlows}
Let $P_i, P_j \in \mathcal{F}_{x,t_i,t_j}[\del_x][\lambda]$ have $\lambda$-degree $\geq 1$, such that Theorem \ref{Isospectral} is satisfied by $\del_{t_i} - P_i$ and $\del_{t_j} - P_j$. On the space of sections over $[\gamma_{\lambda}]$, let $_{[\gamma_{\lambda}](t_i,t_j)}\nabla_{t_i} \cdot = \del_{t_i} + \matr{P}_i$ and $_{[\gamma_{\lambda}](t_i,t_j)}\nabla_{t_j} \cdot = \del_{t_j} + \matr{P}_j$. Then
\[ [ \, [ \del_{t_i} - P_i, \del_{t_j} - P_j ] \, ]_{[\gamma_{\lambda}](t_i,t_j)} = [0]_{[\gamma_{\lambda}](t_i,t_j)} \ . \]
\end{prop}

\prf By Theorem \ref{RankPreserving}, $_{[\gamma_{\lambda}](t_i,t_j)}\nabla_{t_i}$ and $_{[\gamma_{\lambda}](t_i,t_j)}\nabla_{t_j}$ commute with $_{[\gamma_{\lambda}](t_i,t_j)}\nabla_x$, so their commutator commutes with it, as well. By inspection, the first column of the commutator of $_{[\gamma_{\lambda}](t_i,t_j)}\nabla_{t_i}$ with $_{[\gamma_{\lambda}](t_i,t_j)}\nabla_{t_j}$ is $[ \pdiff{P_j}{t_i} - \pdiff{P_i}{t_j} - [P_i,P_j] \, ]_{[\gamma_{\lambda}](t_i,t_j)}$, so this is an element of $\mathcal{E}[\gamma_{\lambda}](t_i,t_j)$. However, it is an element which is polynomial-in-$\lambda$, and by inspection of the basis of $\mathcal{E}[\gamma_\lambda]$, this is only possible if it is identically zero. \quad $\blacksquare$


\subsection{Spectral Deformations of the Kernel of the Tangent Vector Mapping}\label{SpectralEigenring}


In this section, we finish the proof of Theorem \ref{Isospectral} by describing a canonical basis of $\mathcal{E}[\gamma_{\lambda}]$. Without loss of generality, we may ignore the $\llq$time" parameter $t$ in the following.
\begin{thm}\label{CanonicalBasis}
For any spectral-parameterized family $[\gamma_{\lambda}] \in \mathcal{M}$,
\[ \mathcal{E}[\gamma_{\lambda}] = < [1]_{[\gamma_{\lambda}]}, [\Omega]_{[\gamma_{\lambda}]}, \dots, [\Omega^{n-1}]_{[\gamma_{\lambda}]} >_{\R^\lambda} \ , \]
where the series of operators $[\Omega]_{[\gamma_{\lambda}]} = \sum_{j = 0}^{\infty} [\Omega_j]_{[\gamma_{\lambda}]} \lambda^{-j}$ is uniquely determined by the condition that $\Omega_0 = \del_x + \tfrac{1}{n}u_{n-1}$. Moreover, $[\Omega^n]_{[\gamma_{\lambda}]} = [\lambda]_{[\gamma_{\lambda}]}$.
\end{thm}

We note that the theorem stands on its own, as an observation on certain deformations of linear ODEs.
We now make some remarks, whose separate consideration will streamline the proof of Theorem \ref{CanonicalBasis}.

\begin{rmrk}\label{FirstRemark}\upshape

Since $\lambda$ commutes with $\del_x$, it suffices to solve the defining equation $L_{[\gamma_{\lambda}]} \cdot [X]_{[\gamma_{\lambda}]} = [0]_{[\gamma_{\lambda}]}$ for $[X]_{[\gamma_{\lambda}]} = \sum_{j = 0}^{\infty} [X_j]_{[\gamma_{\lambda}]} \lambda^{-j}$. Without loss of generality, we may identify $[X_j]_{[\gamma_{\lambda}]}$ with its representative $X_j$ of minimal degree, which has the expression $X_j = \sum_{k=0}^{n-1} (x_k)_j \del_x^k$. Then by Corollary \ref{ScalarOperatorDuality},
\[ (L_{[\gamma]}-\lambda) X = \hat{X}(L_{[\gamma]} - \lambda) \quad \Leftrightarrow \quad L_{[\gamma]}X = \hat{X}L_{[\gamma]} + \lambda (X-\hat{X}) \ , \]
where $\hat{X} = \Phi_{[\gamma_{\lambda}]}(Y)^*$, for $[X]_{[\gamma_\lambda]} \in \mathcal{D}_{[\gamma_\lambda]}$ and $Y \in \mathcal{D}_{[\gamma_\lambda]}^{\, *}$ both induced by the same solution $\mathcal{X}(0)$ of $\mathcal{X}_x = [\matr{F} + \lambda\matr{E}_{n,1}, \mathcal{X}]$ (see $\S\,$\ref{AnalyticVariations} and $\S\,$\ref{DualCurves}).

This equation is equivalent to the recursive sequence of equations:
\begin{equation}\label{SequenceX}
\left\{ \begin{array}{ccl}
              0 & = & X_0-\hat{X}_0 \ , \\
L_{[\gamma]} \cdot [X_k]_{[\gamma]} & = & [X_{k+1}-\hat{X}_{k+1}]_{[\gamma]} \ , \quad k \geq 0 \ .
\end{array} \right.
\end{equation}
By applying Proposition \ref{Variation} for $h = X_{k+1} - \hat{X}_{k+1}$, this sequence can be recursively solved. \footnote{In \cite{Horocholyn}, this solution procedure was carried out using (\ref{LagrangeMagri}), and it was also shown that the coefficients of $X_k$ must be differential polynomials in the curvatures $u_0, \dots, u_{n-1}$ of $[\gamma]$.}
\end{rmrk}

\begin{rmrk}\label{SecondRemark}\upshape
Let $X$, $\hat{X} = \Phi_{[\gamma_{\lambda}]}(Y)^*$, and $\mathcal{X} = \mathcal{X}(0)$ be as in Remark \ref{FirstRemark}, and let $\mathcal{X} = \big( \mathcal{X}_{ij} \big)_{i,j=0}^{n-1}$ and $\mathcal{X}_{ij} = \sum_{k=0}^{\infty} (\mathcal{X}_{ij})_k \lambda^{-k}$, so that $X = \sum_{k=0}^{n-1} \mathcal{X}_{0,k}\del_x^k$ and $\mathcal{X}_{0,k} = \sum_{j=0}^{\infty} (x_k)_j \lambda^{-j}$. 

It follows from Proposition \ref{GaugeTranslationCor} (for $\beta=\gamma$), Remark \ref{DualGaugeTranslation}, and Corollary \ref{DualVariation}, that each matrix entry $\mathcal{X}_{ij}$ has the following two dual expressions:
\[ \mathcal{X}_{ij} = \delta_j(\del_x^i \cdot [X]_{[\gamma]}) = \big( H_{n-1-j}(L_{[\gamma_{\lambda}]}) \cdot Y\big)\big([\del_x^i]_{[\gamma_{\lambda}]} \big) \ . \]
In particular, if $Y = \sum_{k=0}^{n-1} Y_k \delta_k$ and $Y_k = \sum_{j=0}^{\infty} (Y_k)_j \lambda^{-j}$, then $Y_k = \mathcal{X}_{k,n-1}$ and $\hat{X} = \sum_{k=0}^{n-1} H_{k}(L_{[\gamma_{\lambda}]}) \cdot Y_{n-1-k}$. Consequently, we see that $Y_0 = \mathcal{X}_{0,n-1}$ and $Y_{k+1} = Y_k^{(1)} - u_{n-1}Y_k + \mathcal{X}_{k,n-2}$. Continuing in this fashion, we then deduce that
\begin{equation}\label{y_j}
Y_k = \mathcal{X}_{0,n-1-k} + A_{1,k}(\mathcal{X}_{0,n-k}) + \cdots + A_{k,k}(\mathcal{X}_{0,n-1}) \ ,
\end{equation}
where each $A_{j,k} \in \mathcal{D}_x$ is of degree $j$, with leading coefficient $\tbinom{k}{j}$.

Since $\mathcal{X}_{0,n-1-k} = X_k$, it follows from expanding $\hat{X} = \sum_{k=0}^{n-1} H_{k}(L_{[\gamma_{\lambda}]}) \cdot Y_{n-1-k}$ that the operator
$\Delta_X := \hat{X} - X$ is of the form $\Delta_X = \sum_{m=0}^{n-2} \Delta_m\del_x^m$, where each coefficient $\Delta_m$ has the form:
\[ \Delta_m = \Delta_{n-1,m}(\mathcal{X}_{0,n-1}) + \cdots + \Delta_{m+1,m}(\mathcal{X}_{0,m+1}) \ , \]
such that $\Delta_{k,m} \in \mathcal{D}_x$ of degree $k-m$. In particular, it follows by inspection that $\Delta_{m+1,m} = n\del_x$, meaning that (omitting irrelevant differential polynomials)
\[ \hat{X} = X + \sum_{j=0}^{n-2} \big( n\mathcal{X}_{0,j+1}^{(1)} + \cdots \big) \del_x^j \ . \]
\end{rmrk}

Theorem \ref{CanonicalBasis} follows essentially from the following (cf. \cite[Lemma 1.8.10]{Dickey}):
\begin{lemma}\label{Step1}
Let $[X]_{[\gamma_{\lambda}]} = \sum_{j=0}^{\infty} [X_j]_{[\gamma_{\lambda}]} \lambda^{-j}$ satisfy (\ref{SequenceX}). If $X_0 = \sum_{k=0}^{n-1} (x_k)_0 \del_x^k$ such that none of the coefficients $(x_k)_0$ are constant, then $[X]_{[\gamma_{\lambda}]} \in \R((\lambda^{-1})) \cdot [1]_{[\gamma_{\lambda}]}$.
\end{lemma}

\prf Assume that none of the $(x_j)_0$ are constant, for $1 \leq j \leq n-1$, and let $k$ be the greatest integer between $1$ and $n-1$ such that $(x_k)_0$ is non-zero. By Remark \ref{SecondRemark}, observe that on the sub-diagonal $\mathcal{X}_{i,k-1+i}$ $(i=0, \dots, n-1)$ containing $(x_{k-1})_0$, the $\lambda^0$-terms of the entries are:
\[ (\mathcal{X}_{i,k-1+i})_0 = \begin{cases} i \mathcal{X}_{0,k}^{(1)} + \mathcal{X}_{0,k-1} \ , & 0 \leq i \leq n-2-k \ , \\
(n-k)\mathcal{X}_{0,k}^{(1)} - u_{n-1}\mathcal{X}_{0,k} + \mathcal{X}_{0,k-1} \ , & i = n-1-k \ .
\end{cases} \]
Then $n-1-k$ is the least integer for which $(Y_j)_0$ is non-zero, so the $\lambda^0$ term of $\hat{X}$ is $\hat{X}_0 = \sum_{j=0}^{k} H_j(L_{[\gamma_{\lambda}]}) \cdot (Y_{n-1-j})_0$. In $\hat{X}_0$, the coefficient of $\del_x^k$ is $(Y_{n-1-k})_0= (x_k)_0$, by (\ref{y_j}), and the coefficient of $\del_x^{k-1}$ is $n (x_k)_0^{(1)} + (x_{k-1})_0$, 
by Remark \ref{SecondRemark} and the above observation regarding $(\mathcal{X}_{i,k-1+i})_0$. Hence, $X_0=\hat{X}_0$ implies that $(x_k)_0^{(1)} = 0$, which contradicts the assumption that $k$ was the greatest integer for which $(x_k)_0$ is non-zero, as well as the assumption that $(x_k)_0$ is not constant.  Consequently, if none of the $(x_j)_0$ are constant for all $1 \leq j \leq n-1$, then $[X_0]_{[\gamma_{\lambda}]} = (x_0)_0 [1]_{[\gamma_{\lambda}]}$.

From the form of the matrix $\mathcal{X}$, it can be checked that $\matr{tr}\mathcal{X}$ is of the form $n (x_0)_j$ plus linear differential polynomials in the other $(x_k)_j$, for each $j$. It is a general fact that $\matr{tr}\mathcal{X}$ is constant, so it follows from the above that $(x_0)_0$ is constant. But if $[X_0]_{[\gamma_{\lambda}]} = (x_0)_0 [1]_{[\gamma_{\lambda}]}$ is constant, then by (\ref{SequenceX}):
\[ [0]_{[\gamma]} = L_{[\gamma]} \cdot [X_0]_{[\gamma]} = [X_1 - \hat{X}_1]_{[\gamma]} \quad \Rightarrow \quad 0 = X_1-\hat{X}_1 \ . \]
Consequently, the above arguments apply to the coefficients of $X_1-\hat{X}_1$, and by induction, it follows that $[X]_{[\gamma_{\lambda}]}$ is a constant multiple of the identity. \quad $\blacksquare$ \\

In conclusion, the lemma implies that every $[X]_{[\gamma_{\lambda}]}$ satisfying (\ref{SequenceX}) is determined by the constants in its leading term $[X_0]_{[\gamma_{\lambda}]}$. It follows that there is only one $[\Omega]_{[\gamma_{\lambda}]} = \sum_{j = 0}^{\infty} [\Omega_j]_{[\gamma_{\lambda}]} \lambda^{-j}$ in $\mathcal{E}[\gamma_{\lambda}]$ such that $\Omega_0$ is monic, of degree 1, and has no constant terms besides the leading $\del_x$. (A simple calculation using the trace of the associated matrix $\bar{\Omega}$ then shows that $\Omega_0 = \del_x + \tfrac{1}{n}u_{n-1}$.) 

Since $\mathcal{E}[\gamma_{\lambda}]$ is an associative algebra, it follows that $[\Omega^k]_{[\gamma_{\lambda}]}$, for $1 \leq k \leq n-1$, are the only elements whose $\lambda^0$-term is a monic operator of degree $k$, with no constant terms besides the leading $\del_x^k$.

It remains to show the last identity in Theorem \ref{CanonicalBasis}. Since 
\[ [\Omega^n]_{[\gamma_{\lambda}]} = [\Omega^{n-1}]_{[\gamma_{\lambda}]} \cdot [\Omega]_{[\gamma_{\lambda}]} \ , \]
and the only constant terms in the $\lambda^0$-term of the right-hand side are $\del_x^{n-1}$ and $\del_x$, the only constant term in the leading $\lambda^0$ term of the left-hand side must be $\del_x^n$. But $\del_x^n \equiv \lambda - \sum_{k=0}^{n-1}u_k\del_x^k \mod L_{[\gamma_\lambda]}$, so the only constant term is, in fact, $\lambda$. By Lemma \ref{Step1}, the assertion $[\Omega^n]_{[\gamma_{\lambda}]} = [\lambda]_{[\gamma_{\lambda}]}$ follows.

\begin{cor}
There is a $\matr{V} \in \mathcal{C}^{\infty}(I, \matr{GL}_n\R)$ such that $\bar{\Omega} = \matr{V}(\matr{F}_{vac} + \matr{E}_{n,1}\lambda)\matr{V}^{-1}$, where $\matr{F}_{vac}$ is from the Frenet-Serret equation (\ref{FrenetSerret}) for the $\llq$vacuum" centro-affine curve defined by $L_{vac} = \del_x^n$.
\end{cor}


\end{document}